# Duality for arithmetic Dijkgraaf-Witten theory

Jaro Nicolas Eichler

**Abstract.** Naidu classified pairs of finite groups and 3-cocycles that lead to equivalent Dijkraaf-Witten theories for 3-manifolds. We establish analogous equivalences for arithmetic Dijkgraaf-Witten theory over totally imaginary number fields $F$ containing $n$-th roots of unity, where $n$ is invertible on $X \subseteq \text{spec } \mathcal{O}_F$. For the full ring of integers $X = \text{spec } \mathcal{O}_F$, we give examples with quadratic fields and the quaternion group $Q_8$ where these equivalences fail, but also identify sufficient conditions under which they still hold.

## Contents

## 1. Introduction

Consider a natural number $n$ and a totally imaginary number field $F$ containing the $n$-th roots of unity $\mu_n$. Fix a trivialization $\zeta_n : \mathbb{Z}/n \xrightarrow{\sim} \mu_n$. For a finite set of finite primes $S$, the subscheme

$$X = \text{spec } \mathcal{O}_F \setminus S$$

can be seen as an oriented cobordism from the empty set to the disjoint union of the spectra of local fields

$$Y = \bigsqcup_{v \in S} \text{spec } F_v.$$

For $S$ empty or $n$ invertible on $X$, Kim defines an arithmetic analogue of Dijkgraaf-Witten theory for $X$ [1]. Consider a finite group $G$ and a 3-cocycle $\omega \in Z^3(G, \mathbb{Z}/n)$. For $S$ empty, the arithmetic Dijkgraaf-Witten invariant is defined as the complex number

$$Z^\omega(X) = \sum_{\rho \in \text{hom}(\pi_1 X, G)} \exp\bigl(2\pi i \, \text{inv}_X(\zeta_{n,*}\rho^*\omega)\bigr) \frac{1}{\#G}$$

where $\text{inv}_X : H^3(\pi_1 X, \mu_n) \to \mathbb{Q}/\mathbb{Z}$ is induced by the local invariant maps. By Hermite-Minkowski, the set $\text{hom}(\pi_1 X, G)$ is finite.

In physics, isomorphisms between TQFTs are commonly called dualities. Using extended TQFTs, by Naidu, there is a classification of dualities between Dijkgraaf-Witten theories for 3-manifolds [2]. This point of view relies on the structure of the cobordism category and is not easily transferable to arithmetic Dijkgraaf-Witten theory. The motivation for this work is the question of whether this classification is still applicable to the arithmetic setting.

Inspired by the topological classification, consider duality data of the following form. Let $n$ be a natural number, $K$ a finite group, $A$ a finite, $n$-torsion $K$-module, and $\gamma \in Z^2(K, A)$ be an inhomogeneous, normalized cocycle. Define $A \rtimes_\gamma K$ as the set $A \times K$ with multiplication $(m, l) \cdot$ $(m', l') = (m + l.m' + \gamma(l, l'), ll')$. Similarly consider the dual $\hat{A} = \hom(A, \mathbb{Z}/n)$ and an inhomogeneous, normalized cocycle $\hat{\gamma} \in Z^2(K, \hat{A})$. Assume there is an inhomogeneous, normalized cochain $e \in C^3(K, \mathbb{Z}/n)$ with $de = \gamma \cup \hat{\gamma}$. Consider the groups

$$G = A \rtimes_\gamma K \text{ and } \hat{G} = \hat{A} \rtimes_{\hat{\gamma}} K$$

with projections $a, k$ and $\hat{a}, \hat{k}$ respectively, and the 3-cocycles

$$\omega = k^*e + a \cup k^*\hat{\gamma} \in Z^3(G, \mathbb{Z}/n) \text{ and } \hat{\omega} = \hat{k}^*e + \hat{k}^*\gamma \cup \hat{a} \in Z^3(\hat{G}, \mathbb{Z}/n)$$

Our main result establishes that this duality data indeed leads to equivalences of arithmetic Dijkgraaf-Witten theories:

**Theorem.** Let $L_1, ..., L_r$ be non-archimedean local fields of characteristic zero and $\zeta_v : \mathbb{Z}/n \xrightarrow{\sim} \mu_n$ an isomorphism over $L_v$ for each $v$. Set $Y_v = \spec L_v$ and choose a base point to define its étale fundamental group $\pi_1 Y_v$. Consider $Y = \bigsqcup Y_v$. The duality data induces an isomorphism

$$\Theta : \mathcal{Z}^\omega(Y) \xrightarrow{\sim} \mathcal{Z}^{\hat{\omega}}(Y),$$

which satisfies the following. Let $F$ be a totally imaginary number field and $\zeta : \mathbb{Z}/n \xrightarrow{\sim} \mu_n$ an isomorphism over $F$. Let $S$ be a finite set of finite places of $F$ containing all divisors of $n$. Suppose there is an isomorphism $\bigsqcup \spec L_v \simeq Y$ compatible with the trivializations $\zeta$ and $\zeta_v$. Set $X = \spec \mathcal{O}_F \setminus S$ and choose a base point to define its étale fundamental group $\pi_1 X$. For each $v$, choose a path between the base point of $X$ and the base point of $Y_v$. Then

$$\Theta(\mathcal{Z}^\omega(X)) = \mathcal{Z}^{\hat{\omega}}(X).$$

The proof relies on Artin-Verdier duality and explicit calculations in the group cohomology of $\pi_1 X$ and $\pi_1 Y$ in section 3.2.

If $X = \spec \mathcal{O}_F$ is closed, then $n$ is not invertible on $X$, so the theorem is not applicable. Suppose $\mathcal{A}$ is a finite, locally constant, $n$-torsion sheaf over $X$. The primary issue is that the cup product

$$\cup : H^i(X, \mathcal{A}) \times H^{3-i}(X, \hom(\mathcal{A}, \mathbb{Z}/n)) \to H^3(X, \mathbb{Z}/n)$$

is no longer a perfect pairing. Denote by $H^1(X, \mathcal{A})^\perp \subseteq H^2(X, \hom(\mathcal{A}, \mathbb{Z}/n))$ the orthogonal complement of $H^1(X, \mathcal{A})$. Consider the duality data $(G, \omega)$ and $(\hat{G}, \hat{\omega})$ from before.

**Theorem.** Let $F$ be a totally imaginary number field and $\zeta : \mathbb{Z}/n \xrightarrow{\sim} \mu_n$ an isomorphism over $F$. Set $X = \spec \mathcal{O}_F$. For every homomorphism $\sigma : \pi_1 X \to K$, assume

$$[\sigma^*\gamma] \notin H^1(X, \sigma^*\hat{A})^\perp \setminus \{0\} \text{ and } [\sigma^*\hat{\gamma}] \notin H^1(X, \sigma^*A)^\perp \setminus \{0\}$$

as well as

$$\frac{\#H^1(X, \sigma^*A)}{\#H^0(X, \sigma^*A)} = \frac{\#H^1(X, \sigma^*\hat{A})}{\#H^0(X, \sigma^*\hat{A})}.$$

Then

$$\mathcal{Z}^\omega(X) = \mathcal{Z}^{\hat{\omega}}(X).$$

In section 4.2 these additional assumptions are examined for some quadratic, totally imaginary number fields. As a specific instance of the duality data considered above, take $G$ to be the quaternion group $Q_8$ as the central extension

$$0 \to \mathbb{Z}/2 \to Q_8 \to \mathbb{Z}/2 \times \mathbb{Z}/2 \to 0$$

and $\omega = 0$. This central extension corresponds to a cocycle $\gamma \in Z^2(\mathbb{Z}/2 \times \mathbb{Z}/2, \mathbb{Z}/2)$. On the other hand, take $\hat{G} = (\mathbb{Z}/2)^3$ and $\hat{\omega} = \hat{k}^*\gamma \cup \hat{a}$ where $\hat{a}$ is the projection onto the first factor and $\hat{k}$ is the projection onto the last two factors. A few remarks:
- Since $\omega$ is trivial, $\mathcal{Z}^\omega(X)$ just counts unramified $Q_8$-torsors over $X$
- For any closed, oriented 3-manifold $M$, the Dijkgraaf-Witten invariants $\mathcal{Z}^\omega(M)$ and $\mathcal{Z}^{\hat{\omega}}(M)$ coincide
- This is not necessarily the case for the ring of integers of a number field. By remark 4.15, for $F = \mathbb{Q}\left(\sqrt{-11 \cdot 59 \cdot 107}\right)$, the invariants are $\mathcal{Z}^\omega(X) = 0.5$ and $\mathcal{Z}^{\hat{\omega}}(X) = 3.5$.

Nevertheless there are still number fields for which $\mathcal{Z}^\omega(X) = \mathcal{Z}^{\hat{\omega}}(X)$ holds. Denote by $\beta$ the Bockstein of the exact sequence

$$0 \to \mathbb{Z}/2 \to \mathbb{Z}/4 \to \mathbb{Z}/2 \to 0.$$

We show in lemma 4.21, if the linking form

$$L_X : H^1(X, \mathbb{Z}/2) \times H^1(X, \mathbb{Z}/2) \to \mathbb{Z}/2$$
$$(x, y) \mapsto x \cup \beta(y)$$

is symmetric and $\omega, \hat{\omega}$ are as above, then $X$ satisfies the assumptions of the second theorem. By lemma 4.24, this holds for $F = \mathbb{Q}\left(\sqrt{-p_1 \cdot \ldots \cdot p_r}\right)$ with primes $p_1, \ldots, p_{r-1} \equiv 1 \bmod 4$ and $p_r \equiv 3 \bmod 4$. In particular, this gives a method to compute the number of unramified $Q_8$ extensions.

This paper is adapted from the author's PhD thesis [3], which contains additional topological context.

## 2. Groupoids

**Definition 2.1.** A groupoid is a small category where every morphism is an isomorphism. It is finite if there are finitely many objects and morphisms.

Consider a groupoid $\mathcal{F}$. Since every morphism is invertible, morphisms induce an equivalence relation on the objects of $\mathcal{F}$.

**Example 2.2.** Suppose $G$ is a group. The category $\bullet //G$ that consists of one object $\bullet$ and morphisms $\hom(\bullet, \bullet) = G$ is a groupoid. It is sometimes also denoted by $BG$. More generally, if $G$ acts from the left on a set $S$, define the groupoid $S//G$ by:
- Objects are elements $s \in S$
- Morphisms $g : s \to t$ are elements $g \in G$ such that $t = g.s$

Its equivalence classes are $S/G$.

For a discrete group $\Lambda$, denote by $\mathrm{tors}_\Lambda$ the category of right $\Lambda$-torsors. Objects are free, transitive $\Lambda$-sets and morphisms are $\Lambda$-equivariant maps.

**Definition 2.3.** Let $\mathcal{F}$ be a groupoid and $\Lambda$ a discrete group. A $\Lambda$-torsor $\mathcal{T}$ over $\mathcal{F}$ is a functor

$$\mathcal{T} : \mathcal{F} \to \mathrm{tors}_\Lambda.$$

A complex line bundle over $\mathcal{F}$ is a functor

$$\mathcal{L} : \mathcal{F} \to \mathrm{vect}_\mathbb{C}$$

whose image lies within the subcategory of 1-dimensional vector spaces.

Suppose $F : \mathcal{F} \to \mathcal{H}$ is a functor of groupoids and $\mathcal{L} : \mathcal{H} \to \mathrm{vect}_\mathbb{C}$ a line bundle. Then the pullback $F^*\mathcal{L}$ of $\mathcal{L}$ along $F$ is the line bundle over $\mathcal{F}$ defined as the composition $\mathcal{L} \circ F$.

Recall how fiber products are constructed for groupoids. For functors of groupoids $F : \mathcal{F} \to \mathcal{H}$ and $E : \mathcal{E} \to \mathcal{H}$, their fiber product $\mathcal{F} \times_{\mathcal{H}} \mathcal{E}$ is given by:

- Objects are tuples $(f, e, \psi : F(f) \to E(e))$ consisting of an object $f$ of $\mathcal{F}$, an object $e$ of $\mathcal{E}$, and a morphism $\psi$ of $\mathcal{H}$
- Morphisms are pairs of morphisms $(\varphi : f \to f', \varepsilon : e \to e')$ such that

$$\begin{array}{ccc} F(f) & \xrightarrow{F(\varphi)} & F(f') \\ \psi \downarrow & & \downarrow \psi' \\ E(e) & \xrightarrow{F(\varepsilon)} & E(e') \end{array}$$

commutes

**Definition 2.4.** Let $F : \mathcal{F} \to \mathcal{H}$ be a functor of groupoids. For an object $h$ of $\mathcal{H}$, the fiber over $h$ is

$$\mathcal{F}(h) = \mathcal{F} \times_{\mathcal{H}} \{h\}.$$

So the fiber of $\mathcal{F}$ over $h$ consists of pairs of an object $f$ in $\mathcal{F}$ and a morphism $\psi : F(f) \to h$ of $\mathcal{H}$.

**Definition 2.5.** Let $\mathcal{F}$ be a groupoid and $\mathcal{L} : \mathcal{F} \to \text{vect}_{\mathbb{C}}$ a line bundle. A section $s$ of $\mathcal{L}$ assigns a vector $s(f) \in \mathcal{L}(f)$ to each object $f \in \mathcal{F}$ such that for any morphism $\varphi : f \to f'$, the morphism $\mathcal{L}(\varphi)$ sends $s(f)$ to $s(f')$. The space of all sections is a complex vector space denoted by $H^0(\mathcal{F}, \mathcal{L})$.

If $S$ is a set with a $G$-action, then $H^0(S//G, \mathcal{L})$ consists of $G$-equivariant sections of the pullback of $\mathcal{L}$ along the projection $S \to S//G$.

**Lemma 2.6.** Let $\mathcal{L} : \mathcal{F} \to \text{vect}_{\mathbb{C}}$ and $\mathcal{L}' : \mathcal{F}' \to \text{vect}_{\mathbb{C}}$ be line bundles over groupoids. There is a canonical isomorphism

$$H^0(\mathcal{F}, \mathcal{L}) \otimes H^0(\mathcal{F}', \mathcal{L}') \xrightarrow{\sim} H^0(\mathcal{F} \times \mathcal{F}', p^*\mathcal{L} \otimes p'^*\mathcal{L}')$$

where $p$ and $p'$ are the projections of $\mathcal{F} \times \mathcal{F}'$.

*Proof.* Any object of $\mathcal{F} \times \mathcal{F}'$ can be written as a tuple $(\rho, \rho')$ where $\rho$ is an object of $\mathcal{F}$ and $\rho'$ is an object of $\mathcal{F}'$. Define

$$K : H^0(\mathcal{F}, \mathcal{L}) \otimes H^0(\mathcal{F}', \mathcal{L}') \to H^0(\mathcal{F} \times \mathcal{F}', p^*\mathcal{L} \otimes p'^*\mathcal{L}')$$

as $K(s \otimes s')(\rho, \rho') = s(\rho) \otimes s'(\rho')$. Since morphisms of $\mathcal{F} \times \mathcal{F}'$ are products of morphisms of $\mathcal{F}$ and $\mathcal{F}'$, sections of the tensor product $p^*\mathcal{L} \otimes p'^*\mathcal{L}'$ decompose into a sum of sections of $\mathcal{L}$ and $\mathcal{L}'$. Therefore $K$ is an isomorphism. $\square$

**Lemma 2.7.** Let $F : \mathcal{F} \to \mathcal{H}$ be a functor of groupoids with finite fibers and $\mathcal{L} : \mathcal{H} \to \text{vect}_{\mathbb{C}}$ a line bundle. There is a linear map $F_* : H^0(\mathcal{F}, F^*\mathcal{L}) \to H^0(\mathcal{H}, \mathcal{L})$ defined as

$$F_*(s)(h) = \sum_{[(f, \psi)]} \mathcal{L}(\psi)(s(f)) \frac{1}{\# \text{aut}(f)}$$

where the sum is over equivalence classes of objects $[(f, \psi)]$ of the fiber $\mathcal{F}(h)$ over $h$.

*Proof.* Suppose $\varphi : f \to f'$ is an isomorphism of $\mathcal{F}(h)$ and take a section $s$ of $F^*\mathcal{L}$. Then

$$\begin{array}{ccc} \mathcal{L}(F(f)) & \xrightarrow{\mathcal{L}(F(\varphi))} & \mathcal{L}(F(f')) \\ & \mathcal{L}(\psi) \searrow \quad \swarrow \mathcal{L}(\psi') & \\ & h & \end{array}$$

commutes and sends $s(f)$ to $s(f')$. Since $f$ and $f'$ are equivalent, $\mathrm{aut}(f)$ and $\mathrm{aut}(f')$ are in bijection. In particular they have the same cardinality. So the sum is independent of the choices of objects of the equivalence classes.

If $\sigma : h \to h'$ is a morphism of $\mathcal{H}$, then

$$\mathcal{F}(h) \to \mathcal{F}(h')$$
$$(f, \psi) \mapsto (f, \sigma \circ \psi)$$

is an equivalence of categories. So $F_* s$ is a section of $\mathcal{L}$. $\square$

## 3. Group cohomology

Here the following cochain complexes are considered. The signs are chosen in accordance with [4, Tag 0FNG].

**Definition 3.1.** For a profinite group $\Phi$ and a discrete $\Phi$-module $\mathcal{A}$, let $C^\bullet(\Phi, \mathcal{A})$ be the complex of inhomogeneous, normalized cochains. Elements in degree $i$ are continuous maps $c : \Phi^i \to \mathcal{A}$ satisfying $c(y_1, ..., y_{k-1}, 1, y_{k+1}, ..., y_i) = 0$ for all $1 \leq k \leq i$. It is defined as 0 in negative degrees. The differential is

$$dc(y_1, ..., y_{i+1}) = y_1.c(y_2, ..., y_{i+1}) + \sum_{k=1}^{i} (-1)^k c(y_1, ..., y_k y_{k+1}, ..., y_{i+1}) + (-1)^{i+1} c(y_1, ..., y_i).$$

There is a unique cup product

$$\cup : H^i(\Phi, \mathcal{M}) \otimes H^j(\Phi, \mathcal{M}') \to H^{i+j}(\Phi, \mathcal{M} \otimes \mathcal{M}')$$

that equals the identity $M \otimes N \to M \otimes N$ in degree 0. On inhomogeneous, normalized cochains it is given as [5, Proposition I.1.4.8]

$$(c \cup c')(y_1, ..., y_{i+j}) = c(y_1, ..., y_i) \otimes (y_1 \cdot ... \cdot y_i).c'(y_{i+1}, ..., y_{i+j}).$$

**Lemma 3.2.** The embedding of the complex of inhomogeneous, normalized cochains into the complex of inhomogeneous cochains is a quasi-isomorphism.

*Proof.* Since

$$\mathop{\mathrm{colim}}_{\mathcal{U} \text{ open}} C^\bullet(\Phi/\mathcal{U}, \mathcal{A}^\mathcal{U}) \simeq C^\bullet(\Phi, \mathcal{A}),$$

the claim reduces to the case of a finite group. For finite groups, see [6, Lemma 6.1 and 6.2]. $\square$

Let $\pi$ be a tuple of $r$ profinite groups $\pi_v$ and $\Delta$ a profinite group with homomorphisms $\pi_v \to \Delta$. For a discrete $\Delta$-module $\mathcal{A}$, denote by $\mathcal{A}_{|\pi}$ the tuple of the pullbacks $\mathcal{A}_{|\pi_v}$ along $\pi_v \to \Delta$ and by

$$C^\bullet(\pi, \mathcal{A}_{|\pi}) = \prod_v C^\bullet(\pi_v, \mathcal{A}_{|\pi_v})$$

the product of the complexes where the differential acts componentwise.

**Definition 3.3.** Let $\mathcal{A}$ be a discrete $\Delta$-module and

$$C^\bullet(\Delta, \pi; \mathcal{A}) = \mathrm{cone}\Big(C^\bullet(\Delta, \mathcal{A}) \xrightarrow{|\pi} C^\bullet(\pi, \mathcal{A}_{|\pi})\Big)[-1]$$

the fiber. Explicitly

$$C^i(\Delta, \pi; \mathcal{A}) = C^i(\Delta, \mathcal{A}) \oplus C^{i-1}(\pi, \mathcal{A}_{|\pi})$$

with differential

$$d = \begin{pmatrix} d & 0 \\ -()_{|\pi} & -d \end{pmatrix}$$

acting from the left on columns.

**Lemma 3.4.** Let $\mathcal{A}$ be a discrete $\Delta$-module. The short exact sequence

$$0 \to C^{\bullet-1}(\pi, \mathcal{A}_{|\pi}) \to C^\bullet(\Delta, \pi; \mathcal{A}) \to C^\bullet(\Delta, \mathcal{A}) \to 0$$

induces the long exact sequence

$$\ldots \to H^i(\Delta, \pi; \mathcal{A}) \to H^i(\Delta, \mathcal{A}) \xrightarrow{-()_{|\pi}} H^i(\pi, \mathcal{A}) \to \ldots$$

*Proof.* The short exact sequence of complexes induces a long exact sequence in cohomology. To see that the connecting homomorphism is $-()_{|\pi}$, consider $\alpha \in Z^i(\Delta, \mathcal{A})$. The projection sends $(\alpha, 0) \in C^i(\Delta, \pi; \mathcal{A})$ to $\alpha$ and the differential sends $(\alpha, 0)$ to $(d\alpha, -\alpha_{|\pi})$. Since $\alpha$ is a cocycle, $d\alpha = 0$, so $(d\alpha, -\alpha_{|\pi}) = (0, -\alpha_{|\pi})$. Therefore the connecting homomorphism sends the class $[\alpha]$ to the class $[-\alpha_{|\pi}] \in H^i(\pi, \mathcal{A}_{|\pi})$. □

**Lemma 3.5.** Let $\mathcal{A}$, $\mathcal{B}$, and $\mathcal{C}$ be discrete $\Delta$-modules and $\mathcal{A} \otimes \mathcal{B} \to \mathcal{C}$ a pairing. There are pairings

$$\cup : H^i(\Delta, \pi; \mathcal{A}) \otimes H^j(\Delta, \mathcal{B}) \to H^{i+j}(\Delta, \pi; \mathcal{C})$$
$$(\alpha, \beta) \otimes \gamma \mapsto (\alpha \cup \gamma, \beta \cup \gamma_{|\pi})$$

and

$$\cup : H^j(\Delta, \mathcal{B}) \otimes H^i(\Delta, \pi; \mathcal{A}) \to H^{i+j}(\Delta, \pi; \mathcal{C})$$
$$\gamma \otimes (\alpha, \beta) \mapsto (\gamma \cup \alpha, (-1)^j \gamma_{|\pi} \cup \beta).$$

*Proof.* Consider the map defining cup product

$$\begin{array}{c} \mathcal{A} \otimes \mathcal{B} \\ \downarrow \searrow \\ \mathrm{tot}(C^\bullet(\Delta, \mathcal{A}) \otimes C^\bullet(\Delta, \mathcal{B})) \longrightarrow C^\bullet(\Delta, \mathcal{A} \otimes \mathcal{B}) \longrightarrow C^\bullet(\Delta, \mathcal{C}) \end{array}$$

Then

$$\begin{array}{ccc} \mathrm{tot}(C^\bullet(\Delta, \mathcal{A}) \otimes C^\bullet(\Delta, \mathcal{B})) & \longrightarrow & C^\bullet(\Delta, \mathcal{C}) \\ \downarrow & & \downarrow \\ \prod_v \mathrm{tot}(C^\bullet(\pi_v, \mathcal{A}_{|\pi_v}) \otimes C^\bullet(\Delta, \mathcal{B})) & \longrightarrow & \prod_v C^\bullet(\pi_v, \mathcal{C}_{|\pi_v}) \end{array}$$

commutes. By the functoriality of the cone in the category of cochain comples, by [7, Proposition 6.1], this induces a pairing

$$\mathrm{tot}(C^\bullet(\Delta, \pi; \mathcal{A}) \otimes C^\bullet(\Delta, \mathcal{B})) \to C^\bullet(\Delta, \pi; \mathcal{C}).$$

Take cohomogoly to get the desired pairing and compare with the reference for the formulas. Similarly for the factors switched. □

**Lemma 3.6.** The diagram

$$
\begin{array}{ccc}
H^i(\Delta,\pi;\mathcal{A}) \otimes H^j(\Delta,\mathcal{B}) & \xrightarrow{\cup} & H^{i+j}(\Delta,\pi;\mathcal{C}) \\
{\scriptstyle (-1)^{ij}} \downarrow & & \downarrow \wr \\
H^j(\Delta,\mathcal{B}) \otimes H^i(\Delta,\pi;\mathcal{A}) & \xrightarrow{\cup} & H^{i+j}(\Delta,\pi;\mathcal{C})
\end{array}
$$

commutes.

*Proof.* Suppose $(\alpha,\beta) \in H^i(\Delta,\pi;\mathcal{A})$ and $\gamma \in H^j(\Delta,\mathcal{B})$. By the properties of the cup product in group cohomology, you get

$$\alpha \cup \gamma = (-1)^{ij} \gamma \cup \alpha$$

and

$$\beta \cup \gamma_{|\pi} = (-1)^{(i-1)j} \gamma_{|\pi} \cup \beta = (-1)^{ij}(-1)^j \gamma_{|\pi} \cup \beta.$$

Therefore the claim

$$(\alpha,\beta) \cup \gamma = (-1)^{ij} \gamma \cup (\alpha,\beta)$$

follows. □

## 3.1. Field theory

Let $G$ be a finite group. Denote by

$$\mathrm{ad}_g : G \to G$$
$$s \mapsto gsg^{-1}$$

conjugation. As before let $\pi$ be a tuple of $r$ profinite groups $\pi_v$ and $\Delta$ a profinite group with homomorphisms $\pi_v \to \Delta$. Assume further that $\pi_v$ and $\Delta$ have finite corank.

**Definition 3.7.** The groupoid of fields over $\Delta$ is the groupoid

$$\mathcal{F}_\Delta^G = \mathrm{hom}(\Delta, G)//G$$

where $G$ acts by conjugation. Explicitly:
- Objects are group homomorphisms $\tau : \Delta \to G$
- Morphisms $g : \tau \to \tau'$ are elements $g \in G$ such that $\tau' = g.\tau = \mathrm{ad}_g \circ \tau$

The groupoid of fields over $\pi$ is

$$\mathcal{F}_\pi^G = \prod_v \mathrm{hom}(\pi_v, G)//G^r.$$

Write objects of $\mathcal{F}_\pi^G$ as $\rho : \pi \to G$ representing tuples of $r$ homomorphisms $\rho_v : \pi_v \to G$. Then elements $g \in G^r$ act on $\rho$ componentwise. For a $G$-module $A$, $g \in G$, and $c \in C^{i+1}(G,A)$, the group $G$ acts on cochains by $c.g = g^{-1}.\mathrm{ad}_g^* c$.

**Lemma 3.8.** Let $A$ be a $G$-module. The map $h_g : C^{i+1}(G,A) \to C^i(G,A)$ defined as

$$h_g(c)(y) = \sum_{k=0}^{i} (-1)^k c(y_1, ..., y_k, g^{-1}, \mathrm{ad}_g(y_{k+1}), ..., \mathrm{ad}_g(y_i))$$

is a chain homotopy between the action of $g$ and id. So

$$c.g - c = dh_g(c) + h_g(dc).$$

Moreover, if $g' \in G$ is another element,
$$h_{g'g} = \mathrm{ad}_g^* h_{g'} + h_g.$$

*Proof.* See [8, Proposition 8.1]. □

Let $\omega \in Z^3(G, \mathbb{Z}/n)$. Assume
$$\mathrm{tr}_{\pi_v} : H^2(\pi_v, \mathbb{Z}/n) \xrightarrow{\sim} \frac{1}{n}\mathbb{Z}/\mathbb{Z}$$
is an isomorphism for each $v$ and set $\mathrm{tr}_\pi = \sum_v \mathrm{tr}_{\pi_v}$. Further assume $H^3(\pi_v, \mathbb{Z}/n) = 0$.

**Definition 3.9.** The prequantization bundle $\mathcal{T}_\pi^\omega$ is the $\frac{1}{n}\mathbb{Z}/\mathbb{Z}$-torsor over $\mathcal{F}_\pi^G$ defined by:

- For an object $\rho = (\rho_v : \pi_v \to G)$, let
$$\mathcal{T}_\pi^\omega(\rho) = \mathrm{tr}_{\pi,*} \prod_v d^{-1}(\rho_v^* \omega)/B^2(\pi_v, \mathbb{Z}/n) = \left(\bigoplus_v d^{-1}(\rho_v^* \omega)/B^2(\pi_v, \mathbb{Z}/n)\right)/\sim$$
where $(c_v) \sim (c_{v'})$ if and only if $\mathrm{tr}_\pi(c_{v'} - c_v) = 0$. Since $H^3(\pi_v, \mathbb{Z}/n)$ vanishes, each $d^{-1}(\rho_v^* \omega)$ is non-empty. So this is a
$$\mathrm{tr}_{\pi,*} \bigoplus_v H^2(\pi_v, \mathbb{Z}/n) \xrightarrow{\sim} \frac{1}{n}\mathbb{Z}/\mathbb{Z}$$
torsor.

- For a morphism $g : \rho \to g.\rho$ given by a tuple $g = (g_v) \in G^r$, let $h_g = (h_{g_v})$ be the tuple of chain homotopies from lemma 3.8. Then
$$\mathcal{T}_\pi^\omega(g) : \mathcal{T}_\pi^\omega(\rho) \to \mathcal{T}_\pi^\omega(g.\rho)$$
is defined as addition by $\rho^* h_g(\omega)$.

For an element $c \in \mathcal{T}_\pi^\omega(\rho)$, the calculation
$$d(c + \rho^* h_g(\omega)) = \rho^* \omega + \rho^* dh_g(\omega) = \rho^* \omega + (g.\rho)^* \omega - \rho^* \omega - \rho^* h_g(d\omega) = (g.\rho)^* \omega$$
shows that $c + \rho^* h_g(\omega) \in \mathcal{T}_\pi^\omega(g.\rho)$. So $\mathcal{T}_\pi^\omega(g)$ is well-defined. Suppose $g, g' \in G^r$. Then
$$\rho^* h_{g'g}(\omega) = \rho^* \mathrm{ad}_g^* h_{g'}(\omega) + \rho^* h_g(\omega) = (g.\rho)^* h_{g'}(\omega) + \rho^* h_g(\omega)$$
shows that $\mathcal{T}_\pi^\omega$ is functorial.

Let
$$\mathrm{tr}_\Delta : H^3(\Delta, \pi; \mathbb{Z}/n) \xrightarrow{\sim} \frac{1}{n}\mathbb{Z}/\mathbb{Z}$$
be an isomorphism that commutes with $\mathrm{tr}_\pi$ in the long exact sequence from lemma 3.4. Assume $H^3(\Delta, \mathbb{Z}/n) = 0$.

**Definition 3.10.** Let $\tau : \Delta \to G$ be an object of $\mathcal{F}_\Delta^G$. Since $H^3(\Delta, \mathbb{Z}/n)$ vanishes, there is a cochain $z \in C^2(\Delta, \mathbb{Z}/n)$ such that $dz = \tau^* \omega$. The action of $\tau$ is
$$\mathcal{S}_\Delta^\omega(\tau) = z_{|\pi} \in \mathcal{T}_\pi^\omega(\tau_{|\pi}).$$

If you choose another cochain $z' \in d^{-1}(\tau^* \omega)$, then $\zeta = z - z' \in H^2(\Delta, \mathbb{Z}/n)$. By the long exact sequence from lemma 3.4, the composition

$$H^2(\Delta, \mathbb{Z}/n) \to H^2(\pi, \mathbb{Z}/n) \to H^3(\Delta, \pi; \mathbb{Z}/n)$$

is zero. So, using the assumption on the compatibility of $\mathrm{tr}_\pi$ and $\mathrm{tr}_\Delta$,

$$\mathrm{tr}_\pi(\zeta_{|\pi}) = \mathrm{tr}_\Delta(0) = 0.$$

Therefore $\zeta_{|\pi}$ is trivial in $\mathcal{T}_\pi^\omega(\tau_{|\pi})$, which implies that $\mathcal{S}_\Delta^\omega(\tau)$ is independent of the choice $z$.

Define

$$\exp_\pi : H^2(\pi, \mathbb{Z}/n) \xrightarrow{\mathrm{tr}_\pi} \frac{1}{n}\mathbb{Z}/\mathbb{Z} \xrightarrow{\exp(2\pi i -)} \mathbb{C}^\times.$$

Consider the line bundle

$$\mathcal{L}_\pi^\omega = \mathrm{tr}_{\pi,*} \mathcal{T}_\pi^\omega \times_{\frac{1}{n}\mathbb{Z}/\mathbb{Z}} \mathbb{C}$$

over $\mathcal{F}_\pi^\omega$. For an object $\rho$ of $\mathcal{F}_\pi^\omega$, it is defined as

$$\mathcal{L}_\pi^\omega(\rho) = \{(c, x) \mid c \in \mathcal{T}_\pi^\omega(\rho), x \in \mathbb{C}\}/\sim$$

where

$$(c + \alpha, x) \sim (c, \exp_\pi(\alpha) \cdot x).$$

Addition is defined as

$$(c, x) + (c', x') = (c, x) + (c' + (c - c'), \exp_\pi(c' - c) \cdot x') = (c, x + \exp_\pi(c' - c) \cdot x')$$

and the Hermitian product is defined as

$$<(c, x), (c, x')> = x\overline{x'}$$

where $\overline{x}$ is the complex conjugate. Since

$$x \exp_\pi(-\varepsilon) \overline{x' \exp_\pi(-\varepsilon)} = x\overline{x'},$$

this product is well-defined. Define the exponential

$$\exp_\pi : \mathcal{T}_\pi^\omega(\rho) \to \mathcal{L}_\pi^\omega(\rho)$$
$$c \mapsto (c, 1).$$

**Lemma 3.11.** The exponential of the action $\exp_\pi(\mathcal{S}_\Delta^\omega)$ is a global section of the pullback of $\mathcal{L}_\pi^\omega$ to $\mathcal{F}_\Delta^G$.

*Proof.* If $z \in C^2(\Delta, \mathbb{Z}/n)$ satisfies $dz = \tau^*\omega$, then $z' = z + \tau^*h_g(\omega)$ satisfies

$$dz' = d(z + \tau^*h_g(\omega)) = (g.\tau)^*\omega.$$

So

$$\mathcal{T}_\pi^\omega(g)(z_{|\pi}) = z_{|\pi} + \tau_{|\pi}^* h_g(\omega) = z'_{|\pi}$$

shows that $\mathcal{S}_\Delta^\omega$ is a $G$-equivariant section of the pullback of $\mathcal{T}_\pi^\omega$ to $\mathrm{hom}(\Delta, G)$. Therefore the exponential is a global section. $\square$

**Definition 3.12.** The space of quantum states of $\pi$ is the complex vector space of sections

$$\mathcal{Z}^\omega(\pi) = H^0(\mathcal{F}_\pi^G, \mathcal{L}_\pi^\omega)$$

with Hermitian product

$$< s, s' > = \sum_{[\rho]} < s(\rho), s'(\rho) > \frac{1}{\#\operatorname{aut}(\rho)}$$

where the sum is over equivalence classes $[\rho]$ of $\mathcal{F}_\pi^G$.

**Lemma 3.13.** Let $1 \leq s \leq r$. There is a canonical isomorphism

$$\mathcal{Z}^\omega(\pi_1, ..., \pi_s) \otimes \mathcal{Z}^\omega(\pi_{s+1}, ..., \pi_r) \xrightarrow{\sim} \mathcal{Z}^\omega(\pi).$$

*Proof.* Denote the first tuple by $\pi'$, the second one by $\pi''$, and the projections of $\mathcal{F}_{\pi'}^G \times \mathcal{F}_{\pi''}^G \simeq \mathcal{F}_\pi^G$ by $p'$ and $p''$. Then

$$p'^* \mathcal{L}_{\pi'}^\omega \otimes p''^* \mathcal{L}_{\pi''}^\omega \to \mathcal{L}_\pi^\omega$$
$$((c_{\rho'}), x_{\rho'}) \otimes ((c_{\rho''}), x_{\rho''}) \mapsto ((c_\rho), x_{\rho'} x_{\rho''})$$

defines an isomorphism of line bundles where $(c_\rho) = (c_{\rho'}, c_{\rho''})$ is the concatenated tuple. By lemma 2.6,

$$H^0\bigl(\mathcal{F}_{\pi'}^G, \mathcal{L}_{\pi'}^\omega\bigr) \otimes H^0\bigl(\mathcal{F}_{\pi''}^G, \mathcal{L}_{\pi''}^\omega\bigr) \simeq H^0\bigl(\mathcal{F}_{\pi'}^G \times \mathcal{F}_{\pi''}^G, p'^* \mathcal{L}_{\pi'}^\omega \otimes p''^* \mathcal{L}_{\pi''}^\omega\bigr).$$

This shows the claim. □

By lemma 3.11, the action $\mathcal{S}_\Delta^\omega$ is a section of the pullback of $\mathcal{T}_\pi^\omega$ to $\mathcal{F}_\Delta^G$. By lemma 2.7 the pushforward of $\exp_\pi(\mathcal{S}_\Delta^\omega)$ along $\mathcal{F}_\Delta^G \to \mathcal{F}_\pi^G$ is a section of $\mathcal{L}_\pi^\omega$. Therefore the following definition is well-defined.

**Definition 3.14.** The path integral of $\Delta$ is defined as the section

$$\mathcal{Z}^\omega(\Delta) = \left(\rho \mapsto \sum_{[(\tau, g)]} \mathcal{L}_\pi^G(g)\bigl(\exp_\pi(\mathcal{S}_\Delta^\omega(\tau))\bigr)\right) \in \mathcal{Z}^\omega(\pi)$$

where the sum is over equivalence classes of objects $[(\tau, g)]$ of the fiber $\mathcal{F}_\Delta^G(\rho)$ over $\rho$.

You would expect a factor $1/\#\operatorname{aut}(\tau)$ in the definition of $\mathcal{Z}^\omega(\Delta)$. It vanishes since there is no nontrivial automorphisms in the fibers. If you expand the definition of the summand, you get

$$\mathcal{L}_\pi^G(g)\bigl(\exp_\pi(\mathcal{S}_\Delta^\omega(\tau))\bigr) = \bigl(\mathcal{S}_\Delta^\omega(\tau) + \tau_{|\pi}^* h_g(\omega), 1\bigr)$$

as an element of $\mathcal{L}_\pi^G(\rho)$.

**Lemma 3.15.** Let $H$ be a finite group, $f : H \to G$ a homomorphism, and $f : \mathcal{F}_\pi^H \to \mathcal{F}_\pi^G$ the composition with $f$. Then

$$f^* \mathcal{L}_\pi^\omega = \mathcal{L}_\pi^{f^* \omega}.$$

*Proof.* It is enough to show $f^* \mathcal{T}_\pi^\omega = \mathcal{T}_\pi^{f^* \omega}$ for the prequantization bundle. For an object $\rho$ of $\mathcal{F}_\pi^H$, the sets

$$f^* \mathcal{T}_\pi^\omega(\rho) = \mathcal{T}_\pi^\omega(f \circ \rho) = \operatorname{tr}_{\pi, *} d^{-1}\bigl((f \circ \rho)^* \omega\bigr) / B^2(\pi, \mathbb{Z}/n) = \mathcal{T}_\pi^{f^* \omega}(\rho)$$

coincide. It remains to show that the $H$-action is the same. Take $l \in H$, then $h_l(f^* \omega) = f^* h_{f(l)}(\omega)$. Therefore

$$\rho^* h_l(f^* \omega) = (f \circ \rho)^* h_{f(l)}(\omega)$$

shows the claim. □

**Lemma 3.16.** Let $\omega, \omega' \in Z^3(G, \mathbb{Z}/n)$. There is a canonical isomorphism of line bundles
$$\mathcal{L}_\pi^\omega \otimes \mathcal{L}_\pi^{\omega'} \xrightarrow{\sim} \mathcal{L}_\pi^{\omega+\omega'}.$$

*Proof.* Suppose $\rho : \pi \to G$ is an object of $\mathcal{F}_\pi^G$. On cochains, the map
$$d^{-1}(\rho^*\omega) \oplus d^{-1}(\rho^*\omega') \to d^{-1}(\rho^*(\omega + \omega'))$$
$$(a, b) \mapsto a + b$$

is well-defined and surjective. Since $h_g$ is additive, it induces a morphism of line bundles
$$\mathcal{L}_\pi^\omega \otimes \mathcal{L}_\pi^{\omega'} \to \mathcal{L}_\pi^{\omega+\omega'},$$
which is surjective in each fiber. Therefore this defines an isomorphism of line bundles. □

**Lemma 3.17.** Let $b \in C^2(G, \mathbb{Z}/n)$. There is a canonical trivialization
$$\mathcal{L}_\pi^{db} \simeq \mathbb{C}_{\mathcal{F}_\pi^G}$$
where $\mathbb{C}_{\mathcal{F}_\pi^G}$ denotes the constant line bundle.

*Proof.* It is enough to show that $\rho \mapsto \exp_\pi(\rho^*b)$ is a global section of $\mathcal{L}_\pi^{db}$. This is equivalent to the assignment $\rho_v \mapsto \rho_v^*b$ being $G$-equivariant over $\hom(\pi_v, G)$. Take $g \in G$. Then
$$h_g(db) + dh_g(b) = b.g - b$$
is equivalent to
$$b.g = h_g(db) + dh_g(b) + b.$$

Therefore
$$(g.\rho_v)^*b = \rho_v^*(b.g) = \rho_v^*h_g(db) + \rho_v^*dh_g(b) + \rho_v^*b = \rho_v^*b + \rho_v^*h_g(db) \mod B^2(\pi_v, \mathbb{Z}/n)$$
shows the claim. □

**Corollary 3.18.** Let $b \in C^2(G, \mathbb{Z}/n)$. There is a canonical isomorphism
$$\mathcal{Z}^\omega(\pi) \xrightarrow{\sim} \mathcal{Z}^{\omega+db}(\pi),$$
which sends $\mathcal{Z}^\omega(\Delta)$ to $\mathcal{Z}^{\omega+db}(\Delta)$.

*Proof.* By lemma 3.16 and lemma 3.17, addition by $\rho^*b$ induces an isomorphism of the underlying line bundles. It remains to show that this isomorphism sends $\mathcal{Z}^\omega(\Delta)$ to $\mathcal{Z}^{\omega+db}(\Delta)$. Again it is enough to show this on the level of cochains. If $\tau : \Delta \to G$ is a homomorphism and $z \in C^2(\Delta, \mathbb{Z}/n)$ a cochain such that $dz = \tau^*\omega$, then $d(z + \tau^*b) = \tau^*(\omega + db)$. This shows that the action with respect to $\omega$ is sent to the action with respect to $\omega + db$. □

## 3.2. Duality

Let $K$ be a finite group, $A$ a finite $K$-module, and $\gamma \in Z^2(K, A)$. Define $A \rtimes_\gamma K$ by:
- The underlying set is $A \times K$
- For $(m, l), (m', l') \in A \rtimes_\gamma K$, the multiplication is
$$(m, l) \cdot (m', l') = (m + l.m' + \gamma(l, l'), ll')$$

Denote by $k : A \rtimes_\gamma K \to K$ the projection onto $K$ and by $a : A \rtimes_\gamma K \to A$ the projection onto $A$, which is only a map of sets. The subgroup $A$ acts on $A \rtimes_\gamma K$ non-trivially by conjugation as

$$(m, 1) \cdot (j, l) \cdot (-m, 1) = (m, 1) \cdot (j - l.m, l) = (m + j - l.m, l).$$

**Lemma 3.19.** The sequence

$$1 \to A \to A \rtimes_\gamma K \to K \to 1$$

represents $\gamma$. Further $da = -k^*\gamma$ where $a$ is viewed as a 1-cochain with values in $A$.

*Proof.* The first claim follows from the construction. Take $(m, l), (m', l') \in A \rtimes_\gamma K$. Then

$$\begin{aligned}
da((m, l), (m', l')) &= (m, l).a(m', l') - a((m, l)(m', l')) + a((m, l)) \\
&= l.m' - (m + l.m' + \gamma(l, l')) + m \\
&= -\gamma(l, l')
\end{aligned}$$

shows the second claim. $\square$

Now it is easy to see that $\omega$ and $\hat{\omega}$ are cocycles:

$$d\omega = dk^*e + d(a \cup k^*\hat{\gamma}) = k^*(\gamma \cup \hat{\gamma}) - k^*\gamma \cup k^*\hat{\gamma} - a \cup 0 = 0$$

$$d\hat{\omega} = d\hat{k}^*e + d(\hat{k}^*\gamma \cup \hat{a}) = \hat{k}^*(\gamma \cup \hat{\gamma}) - 0 \cup \hat{a} - \hat{k}^*\gamma \cup \hat{k}^*\hat{\gamma} = 0$$

Suppose $\Phi$ is a profinite group and $\mathcal{A}$ is a discrete $\Phi$-module. Denote by (if it is defined)

$$\chi(\Phi, \mathcal{A}) = \prod_i \#H^i(\Phi, \mathcal{A})^{(-1)^i}$$

the multiplicative Euler characteristic.

As before let $\pi$ be a tuple of $r$ profinite groups $\pi_v$ and $\Delta$ a profinite group with homomorphisms $\pi_v \to \Delta$. View $A$ with the trivial $\Delta$-action.

**Definition 3.20.** If $\chi(\Delta, A)$ is defined, let $\mathcal{Z}^A$ be defined by:
- $\mathcal{Z}^A(\pi) = \mathbb{C}$
- $\mathcal{Z}^A(\Delta) = \chi(\Delta, A)$

Let $\hat{A} = \hom(A, \mathbb{Z}/n)$.

**Lemma 3.21.** For every integer $i$, assume

$$H^i(\Delta, A) \times H^{3-i}(\Delta, \pi; \hat{A}) \overset{\cup}{\to} H^3(\Delta, \pi; \mathbb{Z}/n)$$

is a perfect pairing of finite abelian groups. Then

$$\chi(\Delta, A) = \sqrt{\chi(\pi, A)}.$$

*Proof.* Define

$$\chi(\Delta, \pi; A) = \prod \#H^i(\Delta, \pi; A)^{(-1)^i}.$$

By the perfect pairing

$$H^i(\Delta, A) \times H^{3-i}(\Delta, \pi; \hat{A}) \overset{\cup}{\to} H^3(\Delta, \pi; \mathbb{Z}/n) \simeq \frac{1}{n}\mathbb{Z}/\mathbb{Z}$$

of finite abelian groups, the Euler characteristic is defined and satisfies

$$\chi(\Delta, \pi; \hat{A}) = \chi(\Delta, A)^{-1}.$$

The long exact sequence from lemma 3.4 shows

$$\chi(\Delta, A) = \chi(\Delta, \pi; A) \cdot \chi(\pi, A).$$

Since $\hat{A}$ is isomorphic to $A$ as abelian groups, $\chi(\Delta, \pi; \hat{A}) = \chi(\Delta, \pi; A)$. Therefore $\chi(\Delta, A)^2 = \chi(\pi, A)$ follows. $\square$

**Corollary 3.22.** For every integer $i$, assume

$$H^i(\Delta, A) \times H^{3-i}(\Delta, \pi; \hat{A}) \xrightarrow{\cup} H^3(\Delta, \pi; \mathbb{Z}/n)$$

is a perfect pairing of finite abelian groups. Then multiplication by $\sqrt{\chi(\pi, A)}$ defines an isomorphism

$$\mathbb{C} \xrightarrow{\sim} \mathcal{Z}^A(\pi),$$

which sends 1 to $\mathcal{Z}^A(\Delta)$.

*Proof.* Directly follows from lemma 3.21. $\square$

For an $n$-torsion module $\mathcal{A}$, denote its $\mathbb{Z}/n$-dual by $\hat{\mathcal{A}} = \hom(\mathcal{A}, \mathbb{Z}/n)$.

**Theorem 3.23.** Let $n > 1$ be a natural number and
- $K$ a finite group, $A$ a finite, $n$-torsion $K$-module, and $\hat{A} = \hom(A, \mathbb{Z}/n)$
- $\gamma \in Z^2(K, A), \hat{\gamma} \in Z^2(K, \hat{A})$ inhomogeneous, normalized cocycles, and $e \in C^3(K, \mathbb{Z}/n)$ an inhomogeneous, normalized cochain such that $de = \gamma \cup \hat{\gamma}$
- $G = A \rtimes_\gamma K$ with projections $a$ and $k$ and $\hat{G} = \hat{A} \rtimes_{\hat{\gamma}} K$ with projections $\hat{a}$ and $\hat{k}$
- $\omega = k^*e + a \cup k^*\hat{\gamma} \in Z^3(G, \mathbb{Z}/n)$ and $\hat{\omega} = \hat{k}^*e + \hat{k}^*\gamma \cup \hat{a} \in Z^3(\hat{G}, \mathbb{Z}/n)$

Let $\pi_1, ..., \pi_r$ be profinite groups of finite corank such that $H^3(\pi_v, \mathbb{Z}/n) = 0$. Consider the tuple $\pi = (\pi_1, ..., \pi_r)$. Let

$$\mathrm{tr}_{\pi_v} : H^2(\pi_v, \mathbb{Z}/n) \xrightarrow{\sim} \frac{1}{n}\mathbb{Z}/\mathbb{Z}$$

be isomorphisms and set $\mathrm{tr}_\pi = \sum_v \mathrm{tr}_{\pi_v}$. For every integer $i$, every $v$, and every finite $n$-torsion $\pi_v$-module $\mathcal{A}$, assume $\chi(\pi_v, \mathcal{A})$ is independent of the $\pi_v$-action on $\mathcal{A}$ and

$$H^i(\pi_v, \mathcal{A}) \times H^{2-i}(\pi_v, \hat{\mathcal{A}}) \xrightarrow{\cup} H^2(\pi_v, \mathbb{Z}/n)$$

is a perfect pairing of finite abelian groups. Then there is an isomorphism

$$\Theta : \mathcal{Z}^A(\pi) \otimes \mathcal{Z}^\omega(\pi) \xrightarrow{\sim} \mathcal{Z}^{\hat{\omega}}(\pi).$$

Let $\Delta$ be a profinite group of finite corank with homomorphisms $\pi_v \to \Delta$ such that $H^3(\Delta, \mathbb{Z}/n) = 0$ and

$$\mathrm{tr}_\Delta : H^3(\Delta, \pi; \mathbb{Z}/n) \xrightarrow{\sim} \frac{1}{n}\mathbb{Z}/\mathbb{Z}$$

an isomorphism that commutes with $\mathrm{tr}_\pi$ in the long exact sequence from lemma 3.4. For every integer $i$ and every finite, $n$-torsion $\Delta$-module $\mathcal{A}$, assume $\chi(\Delta, \mathcal{A})$ is independent of the $\Delta$-action on $\mathcal{A}$ and

$$H^i(\Delta, \mathcal{A}) \times H^{3-i}(\Delta, \pi; \hat{\mathcal{A}}) \xrightarrow{\cup} H^3(\Delta, \pi; \mathbb{Z}/n)$$

is a perfect pairing of finite abelian groups. Then

$$\Theta\big(\mathcal{Z}^A(\Delta) \otimes \mathcal{Z}^\omega(\Delta)\big) = \mathcal{Z}^{\hat{\omega}}(\Delta).$$

*Proof.* Since $\mathcal{Z}^A(\pi) = \mathbb{C}$, by composing with the multiplication

$$\mathcal{Z}^A(\pi) \otimes \mathcal{Z}^\omega(\pi) \xrightarrow{\sim} \mathcal{Z}^\omega(\pi),$$

it is enough to construct an isomorphism

$$\Theta : \mathcal{Z}^\omega(\pi) \xrightarrow{\sim} \mathcal{Z}^{\hat{\omega}}(\pi)$$

that sends $\chi(\Delta, A) \cdot \mathcal{Z}^\omega(\Delta)$ to $\mathcal{Z}^{\hat{\omega}}(\Delta)$. First lemma 3.26 constructs $\Theta$ and shows that it is linear. Then lemma 3.31 shows that $\Theta$ is an isomorphism. Finally lemma 3.34 shows

$$\Theta(\chi(\Delta, A) \cdot \mathcal{Z}^\omega(\Delta)) = \mathcal{Z}^{\hat{\omega}}(\Delta),$$

which finishes the proof. $\square$

For the remaining section, assume the notation from theorem 3.23.

**Lemma 3.24.** Let $f : G \to H$ be a surjective group homomorphism and $\varphi$ an object of $\mathcal{F}_\pi^H$. The embedding

$$\{\rho : \pi \to G \mid f \circ \rho = \varphi\}//(\ker f)^r \to \mathcal{F}_\pi^G(\varphi)$$

is an equivalence.

*Proof.* Objects of $\mathcal{F}_\pi^G(\varphi)$ are tuples $(\rho_v : \pi_v \to G, h_v : f \circ \rho_v \to \varphi)$ indexed by $v$ where $h_v \in H$. Take any lift $g_v \in G$ of $h_v$. Then

$$\begin{array}{ccc}
f \circ \rho_v & \xrightarrow{f(g_v)} & f \circ g_v \cdot \rho_v \\
h_v \downarrow & & \downarrow \text{id} \\
\varphi_v & \xrightarrow{\text{id}} & \varphi_v
\end{array}$$

commutes. So the embedding is essentially surjective. Since both are groupoids, it is enough to show that the induced map on automorphism groups is a bijection. But that is obvious, an element $k_v \in G$ such that

$$\begin{array}{ccc}
f \circ \rho_v & \xrightarrow{f(k_v)} & f \circ \rho_v \\
h_v \downarrow & & \downarrow h_v \\
\varphi_v & \xrightarrow{\text{id}} & \varphi_v
\end{array}$$

commutes is exactly an element of the kernel of $f$. $\square$

**Lemma 3.25.** There is a canonical equivalence of groupoids

$$\mathcal{F}_\pi^{G \times_K \hat{G}}(\hat{\rho}) \simeq \mathcal{F}_\pi^G(\hat{k} \circ \hat{\rho}).$$

*Proof.* The embedding $\mathcal{F}_\pi^{G \times_K \hat{G}}$ into $\mathcal{F}_\pi^G \times_{\mathcal{F}_\pi^K} \mathcal{F}_\pi^{\hat{G}}$ is an equivalence. So the top and bottom square in

$$\begin{array}{ccc}
\mathcal{F}_\pi^{G \times_K \hat{G}}(\hat{\rho}) & \longrightarrow & \{\hat{\rho}\} \\
\downarrow & & \downarrow \\
\mathcal{F}_\pi^{G \times_K \hat{G}} & \longrightarrow & \mathcal{F}_\pi^{\hat{G}} \\
\downarrow & & \downarrow \\
\mathcal{F}_\pi^G & \longrightarrow & \mathcal{F}_\pi^K
\end{array}$$

is cartesian. Therefore the whole diagram is cartesian. Since $\{\hat{\rho}\} \to \mathcal{F}_\pi^K$ factors through $\{\hat{k} \circ \hat{\rho}\}$, the claim follows. □

The first step is to construct the linear map

$$\Theta : \mathcal{Z}^\omega(\pi) \to \mathcal{Z}^{\hat{\omega}}(\pi)$$

in the proof of theorem 3.23. Note that you have to compose it with the multiplication

$$\mathcal{Z}^A(\pi) \otimes \mathcal{Z}^\omega(\pi) \twoheadrightarrow \mathcal{Z}^\omega(\pi)$$

to get the map in the statement of theorem 3.23. For $\hat{\rho} : \pi \to \hat{G}$, choose $c_{\hat{\rho}} \in d^{-1}(\hat{\rho}^*\hat{\omega})$. Suppose $s \in \mathcal{Z}^\omega(\pi)$ is a section and write it as $s(\rho) = (c_\rho, x_\rho)$. Define

$$\Theta(c_\rho, x_\rho)(\hat{\rho}) = \left( c_{\hat{\rho}}, \sum_{[\rho]} x_\rho \exp_\pi(c_\rho - \rho^*a \cup \hat{\rho}^*\hat{a} - c_{\hat{\rho}}) \cdot \frac{1}{\#\operatorname{aut}(\rho)} \right)$$

where the sum is over equivalence classes of objects $[\rho]$ of $\{\rho : \pi \to G \mid k \circ \rho = \hat{k} \circ \hat{\rho}\}//A^r$.

**Lemma 3.26.** *The map $\Theta$ is linear.*

*Proof.* The projections of $G \times_K \hat{G}$ induce

$$\begin{array}{ccc}
& \mathcal{F}_\pi^{G \times_K \hat{G}} & \\
{}_p \swarrow & & \searrow {}_{\hat{p}} \\
\mathcal{F}_\pi^G & & \mathcal{F}_\pi^{\hat{G}}
\end{array}$$

by composition. Projecting to $G$ first and then to $K$ is the same as projecting to $\hat{G}$ and then to $K$. By lemma 3.19, $da = -k^*\gamma$ and $d\hat{a} = -\hat{k}^*\hat{\gamma}$. So over $G \times_K \hat{G}$, by implicitly composing with the pullbacks along the projections, you get

$$\omega - d(a \cup \hat{a}) = \omega - da \cup \hat{a} + a \cup d\hat{a} = \omega + k^*\gamma \cup \hat{a} - a \cup \hat{k}^*\hat{\gamma} = \hat{\omega}.$$

View objects of $\mathcal{F}_\pi^{G \times_K \hat{G}}$ as pairs of homomorphisms $\rho : \pi \to G$ and $\hat{\rho} : \pi \to \hat{G}$ such that $k \circ \rho = \hat{k} \circ \hat{\rho}$. For such a tuple $(\rho, \hat{\rho})$, consider the map

$$p^*\mathcal{T}_\pi^\omega(\rho, \hat{\rho}) \to \hat{p}^*\mathcal{T}_\pi^{\hat{\omega}}(\rho, \hat{\rho})$$
$$c \mapsto c - \rho^*a \cup \hat{\rho}^*\hat{a}$$

Composing with $\exp_\pi$, by lemma 3.15 and corollary 3.18, this induces a morphism

$$\beta : p^*\mathcal{L}_\pi^\omega \to \hat{p}^*\mathcal{L}_\pi^{\hat{\omega}}$$

of line bundles over $\mathcal{F}_\pi^{G \times_K \hat{G}}$. Therefore the composition

$$\begin{array}{ccc}
H^0\bigl(\mathcal{F}_\pi^{G \times_K \hat{G}}, p^*\mathcal{L}_\pi^\omega\bigr) & \xrightarrow{\beta} & H^0\bigl(\mathcal{F}_\pi^{G \times_K \hat{G}}, \hat{p}^*\mathcal{L}_\pi^{\hat{\omega}}\bigr) \\
{}_{p^*}\uparrow & & \downarrow {}_{\hat{p}_*} \\
H^0(\mathcal{F}_\pi^G, \mathcal{L}_\pi^\omega) & \dashrightarrow{\Theta} & H^0\bigl(\mathcal{F}_\pi^{\hat{G}}, \mathcal{L}_\pi^{\hat{\omega}}\bigr)
\end{array}$$

defines a linear map $\Theta$. Here $\hat{p}_*$ is the map defined in lemma 2.7 and $p^*(s)(\rho, \hat{\rho}) = s(\rho)$ is the pullback of the section $s$.

It remains to show that this definition of $\Theta$ coincides with the previous one. Suppose $s \in H^0(\mathcal{F}_\pi^G, \mathcal{L}_\pi^\omega)$ is a section. Write $s(\rho) = (c_\rho, x_\rho)$ with $c_\rho \in d^{-1}(\rho^*\omega)$ and $x_\rho \in \mathbb{C}$ to get

$$\beta \circ p^*(s)(\rho, \hat{\rho}) = \beta(c_\rho, x_\rho) = (c_\rho - \rho^* a \cup \hat{\rho}^* \hat{a}, x_\rho).$$

By lemma 3.25 and lemma 3.24,

$$\mathcal{F}_\pi^{G \times_K \hat{G}}(\hat{\rho}) = \mathcal{F}_\pi^G(\hat{k} \circ \hat{\rho}) = \{\rho : \pi \to G \mid k \circ \rho = \hat{k} \circ \hat{\rho}\}//A^r.$$

So, for a section $t \in H^0\left(\mathcal{F}_\pi^{G \times_K \hat{G}}, \hat{p}^* \mathcal{L}_\pi^{\hat{\omega}}\right)$, the pushforward along $\hat{p}$ is

$$\hat{p}_*(t)(\hat{\rho}) = \sum_{[\rho]} t(\rho, \hat{\rho}) \cdot \frac{1}{\#\operatorname{aut}(\rho)}$$

where the sum is over equivalence classes of objects $[\rho]$ of $\{\rho : \pi \to G \mid k \circ \rho = \hat{k} \circ \hat{\rho}\}//A^r$. Choose $c_{\hat{\rho}} \in d^{-1}(\hat{\rho}^* \hat{\omega})$, then

$$\Theta(s)(\hat{\rho}) = \hat{p}_*(\beta \circ p^*(s))(\hat{\rho})$$

$$= \sum_{[\rho]} (c_\rho - \rho^* a \cup \hat{\rho}^* \hat{a}, x_\rho) \cdot \frac{1}{\#\operatorname{aut}(\rho)}$$

$$= \sum_{[\rho]} (c_{\hat{\rho}}, x_\rho \exp_\pi(c_\rho - \rho^* a \cup \hat{\rho}^* \hat{a} - c_{\hat{\rho}})) \cdot \frac{1}{\#\operatorname{aut}(\rho)}$$

$$= \left(c_{\hat{\rho}}, \sum_{[\rho]} x_\rho \exp_\pi(c_\rho - \rho^* a \cup \hat{\rho}^* \hat{a} - x_{\hat{\rho}}) \cdot \frac{1}{\#\operatorname{aut}(\rho)}\right)$$

where the sum is over equivalence classes of objects $[\rho]$ of $\{\rho : \pi \to G \mid k \circ \rho = \hat{k} \circ \hat{\rho}\}//A^r$. $\square$

**Lemma 3.27.** Let $1 \leq s \leq r$. Then

$$\begin{array}{ccc}
\mathcal{Z}^\omega(\pi_1, \ldots, \pi_s) \otimes \mathcal{Z}^\omega(\pi_{s+1}, \ldots, \pi_r) & \xrightarrow{\sim} & \mathcal{Z}^\omega(\pi_1, \ldots, \pi_r) \\
\Theta \otimes \Theta \downarrow & & \Theta \downarrow \\
\mathcal{Z}^{\hat{\omega}}(\pi_1, \ldots, \pi_s) \otimes \mathcal{Z}^{\hat{\omega}}(\pi_{s+1}, \ldots, \pi_r) & \xrightarrow{\sim} & \mathcal{Z}^\omega(\pi_1, \ldots, \pi_r)
\end{array}$$

commutes where the horizontal isomorphisms are defined in lemma 3.13.

*Proof.* The linear map $\Theta$ is the composition $\hat{p}_* \circ \beta \circ p^*$. Both $p^*$ and $\beta$ commute with the isomorphism from lemma 3.13. Since the fields over $\pi$ decompose into a product of the fields over $\pi_1, \ldots, \pi_s$ and $\pi_{s+1}, \ldots, \pi_r$, summing over weighted equivalence classes decomposes in the same way. Therefore the pushforward $\hat{p}_*$ also commutes with the isomorphism. $\square$

**Lemma 3.28.** Let $g = (m, 1) \in G = A \rtimes_\gamma K$. Then $h_g(\omega) = -m \cup k^* \hat{\gamma}$.

*Proof.* Since $e$ is normalized,

$$h_g(k^* e) = k^* h_{k(g)}(e) = k^* h_1(e) = 0.$$

Since $\hat{\gamma}$ is normalized,

$$h_g(a \cup k^* \hat{\gamma})(y, y') = a \cup k^* \hat{\gamma}(g^{-1}, y, y') - a \cup k^* \hat{\gamma}(y, g^{-1}, y') + a \cup k^* \hat{\gamma}(y, y', g^{-1})$$

$$= a(-m) \otimes 1.\hat{\gamma}(k(y), k(y')) - a(y) \otimes k(y).\hat{\gamma}(1, k(y')) + a(y) \otimes k(y).\hat{\gamma}(k(y'), 1)$$

$$= -m \otimes \hat{\gamma}(k(y), k(y')) = -m \cup k^* \hat{\gamma}(y, y').$$

Therefore

$$h_g(\omega) = h_g(k^* e) + h_g(a \cup k^* \hat{\gamma}) = -m \cup k^* \hat{\gamma}$$

shows the claim. $\square$

**Lemma 3.29.** Let $\hat{g} = (\hat{m}, 1) \in \hat{G} = \hat{A} \rtimes_{\hat{\gamma}} K$. Then $h_{\hat{g}}(\hat{\omega}) = -\hat{k}^*\gamma \cup \hat{m}$.

*Proof.* The proof is the same as in lemma 3.28. □

**Lemma 3.30.** Let $\Phi$ be a profinite group and $\mathcal{A}$ a finite $\Phi$-module. Then

$$\frac{\#Z^1(\Phi, \mathcal{A})}{\#\mathcal{A}} = \frac{\#H^1(\Phi, \mathcal{A})}{\#H^0(\Phi, \mathcal{A})}.$$

*Proof.* Examine the defining exact sequence

$$0 \to Z^0(\Phi, \mathcal{A}) \to C^0(\Phi, \mathcal{A}) \to Z^1(\Phi, \mathcal{A}) \to H^1(\Phi, \mathcal{A}) \to 0.$$

Then use $Z^0(\Phi, \mathcal{A}) = H^0(\Phi, \mathcal{A})$ and $C^0(\Phi, \mathcal{A}) = \mathcal{A}$ to get the claim. □

The next step is to show that $\Theta$ in the proof of theorem 3.23 is an isomorphism.

**Lemma 3.31.** For every integer $i$, every $v$, and every finite, $n$-torsion $\pi$-module $\mathcal{A}$, assume $\chi(\pi_v, \mathcal{A})$ is independent of the $\pi_v$-action on $\mathcal{A}$ and

$$H^i(\pi_v, \mathcal{A}) \times H^{2-i}(\pi_v, \hat{\mathcal{A}}) \xrightarrow{\cup} H^2(\pi_v, \mathbb{Z}/n)$$

is a perfect pairing of finite abelian groups. Then $\Theta$ is an isomorphism.

*Proof.* By lemma 3.27, it is enough to show that $\Theta$ is an isomorphism for $\pi = \pi_v$. A good candidate for the inverse is $\hat{\Theta}$ defined by

$$\begin{array}{ccc}
H^0\!\left(\mathcal{F}_\pi^{G \times_K \hat{G}}, p^*\mathcal{L}_\pi^\omega\right) & \xleftarrow{\hat{\beta}} & H^0\!\left(\mathcal{F}_\pi^{G \times_K \hat{G}}, \hat{p}^*\mathcal{L}_\pi^{\hat{\omega}}\right) \\
p_* \downarrow & & \uparrow \hat{p}^* \\
H^0\!\left(\mathcal{F}_\pi^G, \mathcal{L}_\pi^\omega\right) & \xleftarrow{\chi(\pi, \mathcal{A})^{-1} \cdot \hat{\Theta}} & H^0\!\left(\mathcal{F}_\pi^{\hat{G}}, \mathcal{L}_\pi^{\hat{\omega}}\right)
\end{array}$$

where $\hat{\beta}$ is the map induced by $c \mapsto c + \rho^*a \cup \hat{\rho}^*\hat{a}$. Take a section $s \in H^0(\mathcal{F}_\pi^G, \mathcal{L}_\pi^\omega)$ and write $s(\rho) = (c_\rho, x_\rho)$ with $c_\rho \in d^{-1}(\rho^*\omega)$ and $x_\rho \in \mathbb{C}$. Since $s$ is $G$-equivariant,

$$s(g.\rho) = (c_\rho + \rho^*h_g(\omega), x_\rho) = (c_{g.\rho}, x_{g.\rho}).$$

Choose $x_{g.\rho} = x_\rho$ to get $c_{g.\rho} - c_\rho = \rho^*h_g(\omega)$. For $\hat{\rho} \in \mathcal{F}_\pi^{\hat{G}}$, write

$$\Theta(s)(\hat{\rho}) = \left(c_{\hat{\rho}}, \sum_{[\rho]} x_\rho \exp_\pi(c_\rho - \rho^*a \cup \hat{\rho}^*\hat{a} - c_{\hat{\rho}}) \cdot \frac{1}{\#\operatorname{aut}(\rho)}\right) =: (c_{\hat{\rho}}, x_{\hat{\rho}})$$

where the sum is over equivalence classes of objects $[\rho]$ of $\{\rho : \pi \to G \mid k \circ \rho = \hat{k} \circ \hat{\rho}\}/\!/A$. Then

$$\chi(\pi, \mathcal{A})^{-1} \cdot \hat{\Theta}(\Theta(s))(\rho)$$

$$= \left(c_\rho, \sum_{[\hat{\rho}]} x_{\hat{\rho}} \exp_\pi(c_{\hat{\rho}} + \rho^*a \cup \hat{\rho}^*\hat{a} - c_\rho) \cdot \frac{1}{\#\operatorname{aut}(\hat{\rho})}\right)$$

$$= \left(c_\rho, \sum_{[\hat{\rho}]} \sum_{[\rho']} x_{\rho'} \exp_\pi(c_{\rho'} - c_\rho + (\rho^*a - \rho'^*a) \cup \hat{\rho}^*\hat{a}) \cdot \frac{1}{\#\operatorname{aut}(\rho')\#\operatorname{aut}(\hat{\rho})}\right).$$

The section $s$ is trivial if and only if $x_\rho = 0$ for every object $\rho$ of $\mathcal{F}_\pi^G$. If $s$ is trivial, then $\hat{\Theta} \circ \Theta(s)$ is also trivial. Otherwise there is an object $\rho$ of $\mathcal{F}_\pi^G$ such that $s(\rho) \neq 0$. Then $\operatorname{aut}(\rho)$ has to act trivially on $\mathcal{L}_\pi^\omega(\rho)$. An element of $\operatorname{aut}(\rho)$ is a $g \in G$ that is in the centralizer of $\operatorname{im}(\rho)$. It acts on $\mathcal{L}_\pi^\omega(\rho)$ on the first component by addition of the cocycle $\rho^*h_g(\omega)$. This addition is trivial if and only if $\rho^*h_g(\omega)$ is a coboundary. Define $\kappa = k \circ \rho$. For every $(m,1) \in \mathrm{aut}(\rho)$, lemma 3.28 shows

$$\rho^*h_g(\omega) = \rho^*(-m \cup k^*\hat{\gamma}) = -m \cup \kappa^*\hat{\gamma}.$$

So $-m \cup \kappa^*\hat{\gamma}$ is a coboundary. Since

$$H^0(\pi, \kappa^*A) \times H^2(\pi, \kappa^*\hat{A}) \overset{\cup}{\to} H^2(\pi, \mathbb{Z}/n) \overset{\sim}{\to} \frac{1}{n}\mathbb{Z}/\mathbb{Z}$$

is a perfect pairing, the map

$$H^2(\pi, \kappa^*\hat{A}) \overset{\sim}{\to} H^0(\pi, \kappa^*A)^\vee$$
$$[\kappa^*\hat{\gamma}] \mapsto [- \cup \kappa^*\hat{\gamma}]$$

is an isomorphism. So $[\kappa^*\hat{\gamma}] = 0$ if and only if $[m \cup \kappa^*\hat{\gamma}] = 0$ for every $m \in H^0(\pi, \kappa^*A)$. But that is exactly the condition from above with a minus, so $[\kappa^*\hat{\gamma}] = 0$. Therefore the fiber $\mathcal{F}^{\hat{G}}_\pi(\kappa)$ is non-empty. Consider the degree $\#A$ covers $\{\rho' : \pi \to G \mid k \circ \rho' = \kappa\}$ and $\{\hat{\rho} : \pi \to \hat{G} \mid \hat{k} \circ \hat{\rho} = \kappa\}$ of the groupoids that the sums are over. They are $Z^1(\pi, \kappa^*A)$ and $Z^1(\pi, \kappa^*\hat{A})$-torsors respectively. Explicitly an element $\lambda \in Z^1(\pi, \kappa^*A)$ acts as $\lambda\rho'(x) = ((-\lambda + a \circ \rho')(x), \kappa(x))$. Sum over those sets instead and take the degree into account to get

$$\sum_{[\hat{\rho}]} \sum_{[\rho']} x_{\rho'} \exp_\pi(c_{\rho'} - c_\rho + (\rho^*a - \rho'^*a) \cup \hat{\rho}^*\hat{a}) \cdot \frac{1}{\#\mathrm{aut}(\rho)\#\mathrm{aut}(\hat{\rho})}$$

$$= \frac{1}{\#A^2} \sum_{\hat{\rho}} \sum_{\rho'} x_{\rho'} \exp_\pi(c_{\rho'} - c_\rho + (\rho^*a - \rho'^*a) \cup \hat{\rho}^*\hat{a})$$

$$= \frac{1}{\#A^2} \sum_{\lambda} \sum_{\hat{\lambda}} x_{\lambda\rho} \exp_\pi\left(c_{\lambda\rho} - c_\rho + (\rho^*a - (\lambda\rho)^*a) \cup (\hat{\lambda}\hat{\rho})^*\hat{a}\right)$$

$$= \frac{1}{\#A^2} \sum_{\lambda} \sum_{\hat{\lambda}} x_{\lambda\rho} \exp_\pi(c_{\lambda\rho} - c_\rho + \lambda \cup (\hat{\rho}^*\hat{a} - \hat{\lambda}))$$

$$= \frac{1}{\#A^2} \sum_{\lambda} x_{\lambda\rho} \exp_\pi(c_{\lambda\rho} - c_\rho + \lambda \cup \hat{\rho}^*\hat{a}) \sum_{\hat{\lambda}} \exp_\pi(-\lambda \cup \hat{\lambda})$$

where the sums are over $\lambda \in Z^1(\pi, \kappa^*A)$ and $\hat{\lambda} \in Z^1(\pi, \kappa^*\hat{A})$. Since

$$H^1(\pi, \kappa^*A) \times H^1(\pi, \kappa^*\hat{A}) \overset{\cup}{\to} H^2(\pi, \mathbb{Z}/n) \overset{\exp_\pi}{\to} \mathbb{C}^\times$$

is a perfect pairing,

$$\sum_{\hat{\lambda}} \exp_\pi(-\lambda \cup \hat{\lambda})$$

sums over all characters of $H^1(\pi, \kappa^*A)$. So it vanishes if and only if $\lambda \notin B^1(\pi, \kappa^*A)$. Therefore the sum reduces to

$$\frac{\#Z^1(\pi, \kappa^*\hat{A})}{\#A^2} \sum_{\lambda} x_{\lambda\rho} \exp_\pi(c_{\lambda\rho} - c_\rho + \lambda \cup \hat{\rho}^*\hat{a})$$

where the sum is over $\lambda \in B^1(\pi, \kappa^*A)$. If $\lambda$ is a coboundary, it is given by $\lambda(x) = \kappa(x).m - m$ for an $m \in A$ up to elements in $Z^0(\pi, \kappa^*A)$. By lemma 3.28, the element $m$ satisfies $\rho^*h_m(\omega) = -m \cup \kappa^*\hat{\gamma}$. Since

$$d(m \cup \hat{\rho}^*\hat{a}) = dm \cup \hat{\rho}^*\hat{a} - m \cup \kappa^*\hat{\gamma} = (\kappa.m - m) \cup \hat{\rho}^*\hat{a} - m \cup \kappa^*\hat{\gamma}$$

and
$$\lambda\rho = (-\lambda + a \circ \rho, \kappa) = (-\kappa.m + m + a \circ \rho, \kappa) = m.\rho,$$
the above sum is
$$\frac{\#Z^1(\pi, \kappa^*\hat{A})}{\#A^2} \sum_{[m]} x_{m.\rho} \exp_\pi(c_{m.\rho} - c_\rho + (\kappa.m - m) \cup \hat{\rho}^*\hat{a})$$
$$= \frac{\#Z^1(\pi, \kappa^*\hat{A})}{\#A^2} \sum_{[m]} x_\rho \exp_\pi(\rho^* h_m(\omega) + (\kappa.m - m) \cup \hat{\rho}^*\hat{a})$$
$$= \frac{\#Z^1(\pi, \kappa^*\hat{A})}{\#A^2} \sum_{[m]} x_\rho \exp_\pi(-m \cup \kappa^*\hat{\gamma} + (\kappa.m - m) \cup \hat{\rho}^*\hat{a})$$
$$= \frac{\#Z^1(\pi, \kappa^*\hat{A})}{\#A^2} \sum_{[m]} x_\rho \exp_\pi(d(m \cup \hat{\rho}^*\hat{a}))$$
$$= \frac{\#Z^1(\pi, \kappa^*\hat{A})}{\#A \cdot \#Z^0(\pi, \kappa^*A)} x_\rho$$

where the sum is over classes $[m] \in A/Z^0(\pi, \kappa^*A)$. By lemma 3.30 and the perfect pairing, this equals

$$\frac{\#H^1(\pi, \kappa^*\hat{A})}{\#H^0(\pi, \kappa^*\hat{A}) \cdot \#H^0(\pi, \kappa^*A)} x_\rho = \frac{\#H^1(\pi, \kappa^*\hat{A})}{\#H^0(\pi, \kappa^*\hat{A}) \cdot \#H^2(\pi, \kappa^*\hat{A})} x_\rho = \chi(\pi, \kappa^*\hat{A})^{-1} x_\rho.$$

Since $\kappa^*\hat{A} \simeq A$ as abelian groups, by assumption,
$$\chi(\pi, \kappa^*\hat{A}) = \chi(\pi, A).$$

So $\hat{\Theta} \circ \Theta = \text{id}$ and similarly $\Theta \circ \hat{\Theta} = \text{id}$. Therefore $\Theta$ is an isomorphism. □

**Lemma 3.32.** Let $\sigma : \Delta \to K$. Then
$$d^{-1}(\sigma^*\gamma) \to \{\tau : \Delta \to G \mid k \circ \tau = \sigma\}$$
$$b \mapsto (-b, \sigma)$$
is a bijection.

*Proof.* Consider any map of sets
$$\Delta \to A \rtimes_\gamma K$$
$$y \mapsto (b(y), \sigma(y))$$
that continues $\sigma$. It is a homomorphism if and only if
$$b(yy') = b(y) + \sigma(y).b(y') + \gamma(\sigma(y), \sigma(y')),$$
which is exactly the condition for $\sigma^*\gamma = -db$. Therefore
$$d^{-1}(\sigma^*\gamma) \to \{\tau : \Delta \to G \mid k \circ \tau = \sigma\}$$
$$b \mapsto (-b, \sigma)$$
is a bijection. □

Suppose $\mathcal{A}$ is a discrete $\Delta$-module. Recall that elements of $Z^2(\Delta, \pi; \mathcal{A})$ are tuples $(a, b) \in C^2(\Delta, \mathcal{A}) \oplus C^1(\pi, \mathcal{A}_{|\pi})$ such that $da = 0$ and $-a_{|\pi} - db = 0$.

**Lemma 3.33.** Let $\sigma : \Delta \to K$ and $\rho : \pi \to G$ such that $k \circ \rho = \sigma_{|\pi}$. Then

$$d^{-1}(\sigma^*\gamma, \rho^*a) \to \{(\tau : \Delta \to G, m : \tau_{|\pi} \to \rho \mid k \circ \tau = \sigma, m \in A^r\}$$
$$(b, m) \mapsto ((-b, \sigma), m)$$

is a bijection.

*Proof.* Suppose $(b, m) \in C^1(\Delta, \pi; \sigma^*A)$ such that $db = \sigma^*\gamma$. By lemma 3.32, this defines a homomorphism $\tau : \Delta \to G$ as $\tau(x) = (-b(x), \sigma(x))$. Since

$$d(b, m) = (db, -b_{|\pi} - dm) = (db, -b_{|\pi} - \sigma_{|\pi}.m + m) = (db, -b_{|\pi} - (k \circ \rho).m + m),$$

it remains to show that $-b_{|\pi} - (k \circ \rho).m + m = \rho^*a$ if and only if $m.(-b_{|\pi}, \sigma_{|\pi}) = m.\tau_{|\pi} = \rho$. Since $k \circ m.(-b_{|\pi}, \sigma_{|\pi}) = k \circ \rho$, the second condition can be replaced by $a \circ m.(-b_{|\pi}, k \circ \rho) = \rho^*a$. Then

$$a \circ m.(-b_{|\pi}, k \circ \rho) = a \circ (m - b_{|\pi} - (k \circ \rho).m, k \circ \rho) = m - b_{|\pi} - (k \circ \rho).m$$

shows the claim. □

To finish the proof of theorem 3.23, it is left to show that the path integral maps to the path integral.

**Lemma 3.34.** For every integer $i$ and every finite, $n$-torsion $\Delta$-module $\mathcal{A}$, assume $\chi(\Delta, \mathcal{A})$ is independent of the $\Delta$-action on $\mathcal{A}$ and

$$H^i(\Delta, \mathcal{A}) \times H^{3-i}(\Delta, \pi; \hat{A}) \xrightarrow{\cup} H^3(\Delta, \pi; \mathbb{Z}/n)$$

is a perfect pairing of finite abelian groups. Then

$$\Theta(\mathcal{Z}^A(\Delta) \cdot \mathcal{Z}^\omega(\Delta)) = \mathcal{Z}^{\hat{\omega}}(\Delta).$$

*Proof.* Choose $c_\rho \in d^{-1}(\rho^*\omega)$ for every object $\rho$ of $\mathcal{F}_\pi^G$ and $z_\tau \in d^{-1}(\tau^*\omega)$ for every object $\tau$ in $\mathcal{F}_\Delta^G$. Similarly choose $c_{\hat{\rho}}$ and $z_{\hat{\tau}}$. Then

$$\mathcal{Z}^\omega(\Delta)(\rho) = \sum_{[(\tau,g)]} \mathcal{L}_\pi^\omega(g)\left(\exp_\pi\left(\mathcal{S}_\Delta^\omega(\tau)\right)\right)$$

$$= \sum_{[(\tau,g)]} \left(\mathcal{S}_\Delta^\omega(\tau) + \tau_{|\pi}^*h_g(\omega), 1\right)$$

$$= \sum_{[(\tau,g)]} \left(c_\rho, \exp_\pi\left(z_{\tau|\pi} + \tau_{|\pi}^*h_g(\omega) - c_\rho\right)\right)$$

$$= \left(c_\rho, \sum_{[(\tau,g)]} \exp_\pi\left(z_{\tau|\pi} + \tau_{|\pi}^*h_g(\omega) - c_\rho\right)\right) =: (c_\rho, x_\rho)$$

where the sum is over equivalence classes of objects $[(\tau, g)]$ of $\mathcal{F}_\Delta^G(\rho)$. Take an object $\hat{\rho}$ of $\mathcal{F}_\pi^{\hat{G}}$ and denote $\kappa = \hat{k} \circ \hat{\rho}$. Apply $\Theta$ to $\mathcal{Z}^\omega(\Delta)$ to get

$$\Theta(\mathcal{Z}^\omega(\Delta))(\hat{\rho}) = \left(c_{\hat{\rho}}, \sum_{[\rho]} x_\rho \exp_\pi(c_\rho - \rho^*a \cup \hat{\rho}^*\hat{a} - c_{\hat{\rho}}) \frac{1}{\#\operatorname{aut}(\rho)}\right)$$

$$= \left(c_{\hat{\rho}}, \sum_{[\rho]} \sum_{[(\tau,g)]} \exp_\pi\left(z_{\tau|\pi} + \tau_{|\pi}^*h_g(\omega) - \rho^*a \cup \hat{\rho}^*\hat{a} - c_{\hat{\rho}}\right) \frac{1}{\#\operatorname{aut}(\rho)}\right) =: (c_{\hat{\rho}}, x_{\hat{\rho}})$$

where the first sum is over equivalence classes of objects $[\rho]$ of $\{\rho : \pi \to G \mid k \circ \rho = \kappa\}//A^r$ and the second sum is over equivalence classes of objects $[(\tau, g)]$ of $\mathcal{F}^G_\Delta(\rho)$. The claim is that this multiplied by $\chi(\Delta, A)$ equals

$$\mathcal{Z}^{\hat{\omega}}(\Delta)(\hat{\rho}) = \left(c_{\hat{\rho}}, \sum_{[(\hat{\tau}, \hat{g})]} \exp_\pi \left(z_{\hat{\tau}|\pi} + \hat{\tau}^*_{|\pi} h_{\hat{g}}(\hat{\omega}) - c_{\hat{\rho}}\right)\right) =: (c_{\hat{\rho}}, \hat{x}_{\hat{\rho}})$$

where the sum is over equivalence classes of objects $[(\hat{\tau}, \hat{g})]$ of $\mathcal{F}^{\hat{G}}_\Delta(\hat{\rho})$. It is enough to compare $x_{\hat{\rho}}$ to $\hat{x}_{\hat{\rho}}$. Note that $x_{\hat{\rho}}$ and $\hat{x}_{\hat{\rho}}$ both vanish if $\mathcal{F}^K_\Delta(\kappa)$ is empty. Since $\Theta(\mathcal{Z}^\omega(\Delta))$ and $\mathcal{Z}^{\hat{\omega}}(\Delta)$ are sections, it is enough to compare them in any representative $\hat{\rho}$ per equivalence class. Suppose $(\sigma, l)$ is an object of $\mathcal{F}^K_\Delta(\kappa)$. Take any lift $\hat{l} \in \hat{G}^r$ of $l \in K^r$ and replace $\hat{\rho}$ by $\hat{l}^{-1}.\hat{\rho}$. Then $(\sigma, 1)$ is an object of $\mathcal{F}^K_\Delta(\kappa)$.

First examine $x_{\hat{\rho}}$. The set $\text{hom}(\pi, G)(\kappa) = \{\rho : \pi \to G \mid k \circ \rho = \kappa\}$ is a degree $\#A^r$ cover of $\mathcal{F}^G_\pi(\kappa) = \{\rho : \pi \to G \mid k \circ \rho = \kappa\}//A^r$. Sum over elements of this cover instead and compensate for the degree by multiplying with $1/\#A^r$. Since this is a set, summing over $\text{hom}(\pi, G)(\kappa)$ and then over equivalence classes of $\mathcal{F}^G_\Delta(\rho)$ is the same as summing over equivalence classes of the groupoid

$$\mathcal{F}^G_\Delta \times_{\mathcal{F}^G_\pi} \text{hom}(\pi, G)(\kappa) = \{\tau : \Delta \to G, \rho : \pi \to G, m : \tau_{|\pi} \to \rho \mid k \circ \rho = \kappa\}$$

where morphisms $g : \tau \to \tau'$ are elements $g \in A$ such that

$$\begin{array}{ccc} \tau_{|\pi} & \xrightarrow{g} & \tau'_{|\pi} \\ m \downarrow & & \downarrow m' \\ \rho & \xrightarrow{\text{id}} & \rho' \end{array}$$

commutes. Its equivalence classes are fibered over the set of objects of $\mathcal{F}^K_\Delta(\kappa)$. The fiber over $(\sigma, 1)$ is

$$\{\tau : \Delta \to G, m \in A^r \mid k \circ \tau = \sigma\}$$

with the same morphisms as before. Therefore the equivalence classes are in bijection with

$$(\{\tau : \Delta \to G \mid k \circ \tau = \sigma\} \times A^r)/A.$$

By lemma 3.32, this is the same as

$$\left(d^{-1}(\sigma^*\gamma) \times A^r\right)/A.$$

Note that $d^{-1}(\sigma^*\gamma)$ is a $Z^1(\Delta, \sigma^*A)$-torsor, where $\lambda \in Z^1(\Delta, \sigma^*A)$ acts on $\tau : \Delta \to G$ as

$$\lambda\tau = (\tau - \lambda, \sigma).$$

Also note that the $A$-action is free. Therefore, if you sum over the elements in $d^{-1}(\sigma^*\gamma) \times A^r$ instead, you get an extra factor $1/\#A$. Denote by $x_{\hat{\rho}}(\sigma)$ the restriction of the sum $x_{\hat{\rho}}$ to the fiber over $\sigma$. If the cohomology class $[\sigma^*\gamma] \in H^2(\Delta, \sigma^*A)$ does not vanish, the sum is empty. Otherwise take a homomorphism $\tau : \Delta \to G$ such that $k \circ \tau = \sigma$. Since $H^3(\Delta, \mathbb{Z}/n) = 0$, for $\lambda \in Z^1(\Delta, \sigma^*A)$, choose $l_\lambda \in C^2(\Delta, \mathbb{Z}/n)$ such that $dl_\lambda = \lambda \cup \sigma^*\hat{\gamma}$ and define $z_{\lambda\tau} = z_\tau - l_\lambda$. Since

$$dz_{\lambda\tau} = dz_\tau - dl_\lambda = \tau^*\omega - \lambda \cup \sigma^*\hat{\gamma} = \sigma^*e + (\tau^*a - \lambda) \cup \sigma^*\hat{\gamma} = (\lambda\tau)^*\omega,$$

this is well-defined. By lemma 3.28, for $m \in A^r$,

$$(\lambda\tau_{|\pi})^* h_m(\omega) = (\lambda\tau_{|\pi})^*(-m \cup k^*\gamma) = -m \cup \kappa^*\hat{\gamma}$$

and by lemma 3.32, $\lambda\tau = (a \circ \tau - \lambda, \sigma)$. So

$$m.(\lambda\tau_{|\pi}) = m.(a \circ \tau_{|\pi} - \lambda_{|\pi}, \kappa) = (a \circ \tau_{|\pi} + m - \kappa.m - \lambda_{|\pi}, \kappa) = (\tau_{|\pi}^* a - dm - \lambda_{|\pi}, \kappa).$$

Then

$$x_{\hat{\rho}}(\sigma) = \frac{1}{\#A^{r+1}} \mathbb{1}_{[\sigma^*\gamma]=0} \sum_{\lambda,m} \exp_\pi\big(z_{\lambda\tau|\pi} + (\lambda\tau_{|\pi})^* h_m(\omega) - (m.(\lambda\tau_{|\pi}))^* a \cup \hat{\rho}^* \hat{a} - c_{\hat{\rho}}\big)$$

$$= \frac{1}{\#A^{r+1}} \mathbb{1}_{[\sigma^*\gamma]=0} \sum_{\lambda,m} \exp_\pi\big(z_{\tau|\pi} - l_{\lambda|\pi} - m \cup \kappa^*\hat{\gamma} - (\tau_{|\pi}^* a - dm - \lambda_{|\pi}) \cup \hat{\rho}^*\hat{a} - c_{\hat{\rho}}\big)$$

$$= \frac{1}{\#A^{r+1}} \mathbb{1}_{[\sigma^*\gamma]=0} \exp_\pi\big(z_{\tau|\pi} - \tau_{|\pi}^* a \cup \hat{\rho}^* \hat{a} - c_{\hat{\rho}}\big)$$
$$\cdot \sum_m \exp_\pi(dm \cup \hat{\rho}^*\hat{a} - m \cup \kappa^*\hat{\gamma}) \sum_\lambda \exp_\pi(-l_{\lambda|\pi} + \lambda_{|\pi} \cup \hat{\rho}^*\hat{a})$$

$$= \frac{1}{\#A^{r+1}} \mathbb{1}_{[\sigma^*\gamma]=0} \exp_\pi\big(z_{\tau|\pi} - \tau_{|\pi}^* a \cup \hat{\rho}^* \hat{a} - c_{\hat{\rho}}\big)$$
$$\cdot \sum_m \exp_\pi(d(m \cup \hat{\rho}^*\hat{a})) \sum_\lambda \exp_\pi(-l_{\lambda|\pi} + \lambda_{|\pi} \cup \hat{\rho}^*\hat{a})$$

$$= \frac{1}{\#A} \mathbb{1}_{[\sigma^*\gamma]=0} \exp_\pi\big(z_{\tau|\pi} - \tau_{|\pi}^* a \cup \hat{\rho}^* \hat{a} - c_{\hat{\rho}}\big) \sum_\lambda \exp_\pi(-l_{\lambda|\pi} + \lambda_{|\pi} \cup \hat{\rho}^*\hat{a})$$

where the sum is over $\lambda \in Z^1(\Delta, \sigma^*A)$ and $m \in A^r$. Since

$$d(l_\lambda, 0) = (\lambda \cup \sigma^*\hat{\gamma}, -l_{\lambda|\pi}),$$

by lemma 3.5, you get

$$[\lambda] \cup [(\sigma^*\hat{\gamma}, \hat{\rho}^*\hat{a})] = [(\lambda \cup \sigma^*\hat{\gamma}, -\lambda_{|\pi} \cup \hat{\rho}^*\hat{a})] = [(0, -\lambda_{|\pi} \cup \hat{\rho}^*\hat{a} + l_{\lambda|\pi})]$$

in $H^3(\Delta, \pi; \mathbb{Z}/n)$. So, by the compatibility of $\exp_\pi$ and $\exp_\Delta$, and the perfect pairing

$$H^1(\Delta, \sigma^*A) \times H^2(\Delta, \pi; \sigma^*\hat{A}) \to H^3(\Delta, \pi; \sigma^*A),$$

the second sum reduces to

$$\sum_\lambda \exp_\pi(-l_{\lambda|\pi} + \lambda_{|\pi} \cup \hat{\rho}^*\hat{a}) = \sum_\lambda \exp_\Delta(-\lambda \cup (\sigma^*\hat{\gamma}, \hat{\rho}^*\hat{a})) = \mathbb{1}_{[(\sigma^*\hat{\gamma}, \hat{\rho}^*\hat{a})]=0} \#Z^1(\Delta, \sigma^*A).$$

Therefore

$$x_{\hat{\rho}}(\sigma) = \frac{\#Z^1(\Delta, \sigma^*A)}{\#A} \mathbb{1}_{[\sigma^*\gamma]=0} \mathbb{1}_{[(\sigma^*\hat{\gamma}, \hat{\rho}^*\hat{a})]=0} \exp_\pi\big(z_{\tau|\pi} - \tau_{|\pi}^* a \cup \hat{\rho}^* \hat{a} - c_{\hat{\rho}}\big).$$

Now consider $\hat{x}_{\hat{\rho}}$. The fiber of $\mathcal{F}_\Delta^{\hat{G}}(\hat{\rho})$ over $\sigma$ is

$$\{(\hat{\tau} : \Delta \to \hat{G}, \widehat{m} : \hat{\tau}_{|\pi} \to \hat{\rho}) \mid \hat{k} \circ \hat{\tau} = \sigma, \widehat{m} \in \hat{A}^r\}$$

with morphisms $\hat{g} \in \hat{A}$ such that the obvious diagram commutes. By lemma 3.33, its equivalence classes are in bijection with

$$d^{-1}(\sigma^*\hat{\gamma}, \hat{\rho}^*\hat{a})/B^1(\Delta, \pi; \sigma^*\hat{A}).$$

Note that this is an $H^1(\Delta, \pi; \sigma^*\hat{A})$-torsor, where $(\hat{\lambda}, \hat{\mu})$ acts on $(\hat{\tau}, \widehat{m})$ as $(\hat{\lambda}\hat{\tau}, \widehat{m} + \hat{\mu})$. Denote by $\hat{x}_{\hat{\rho}}(\sigma)$ the restriction of the sum $\hat{x}_{\hat{\rho}}$ to the fiber over $\sigma$. If the cohomology class $[(\sigma^*\hat{\gamma}, \hat{\rho}^*\hat{a})] \in H^2(\Delta, \pi; \sigma^*\hat{A})$ does not vanish, the sum is empty. Otherwise take $\hat{\tau} : \Delta \to \hat{G}$ such that $\hat{k} \circ \hat{\tau} = \sigma$ and an element $\hat{m} \in \hat{A}^r$ viewed as a morphism $\hat{m} : \hat{\tau}_{|\pi} \to \hat{\rho}$. As before, for $\hat{\lambda} \in Z^1(\Delta, \sigma^*\hat{A})$, choose $l_{\hat{\lambda}} \in C^2(\Delta, \mathbb{Z}/n)$ such that $dl_{\hat{\lambda}} = \sigma^*\gamma \cup \hat{\lambda}$ and define $z_{\hat{\lambda}\hat{\tau}} = z_{\hat{\tau}} - l_{\hat{\lambda}}$. By lemma 3.29,

$$(\hat{\lambda}\hat{\tau}_{|\pi})^* h_{\hat{m}}(\hat{\omega}) = -\kappa^*\gamma \cup \hat{m}.$$

Then

$$\hat{x}_{\hat{\rho}}(\sigma) = \mathbb{1}_{[(\sigma^*\hat{\gamma}, \hat{\rho}^*\hat{a})]=0} \sum_{[(\hat{\lambda},\hat{\mu})]} \exp_\pi\left(z_{\hat{\lambda}\hat{\tau}|\pi} + (\hat{\lambda}\hat{\tau}_{|\pi})^* h_{\hat{m}+\hat{\mu}}(\hat{\omega}) - c_{\hat{\rho}}\right)$$

$$= \mathbb{1}_{[(\sigma^*\hat{\gamma}, \hat{\rho}^*\hat{a})]=0} \sum_{[(\hat{\lambda},\hat{\mu})]} \exp_\pi\left(z_{\hat{\tau}|\pi} - l_{\hat{\lambda}|\pi} - \kappa^*\gamma \cup (\hat{m} + \hat{\mu}) - c_{\hat{\rho}}\right)$$

$$= \mathbb{1}_{[(\sigma^*\hat{\gamma}, \hat{\rho}^*\hat{a})]=0} \exp_\pi(z_{\hat{\tau}|\pi} - \kappa^*\gamma \cup \hat{m} - c_{\hat{\rho}}) \sum_{[(\hat{\lambda},\hat{\mu})]} \exp_\pi\left(-l_{\hat{\lambda}|\pi} - \kappa^*\gamma \cup \hat{\mu}\right)$$

where the sum is over cohomology classes $[(\hat{\lambda}, \hat{\mu})] \in H^1(\Delta, \pi; \sigma^*\hat{A})$. Since

$$d(l_{\hat{\lambda}}, 0) = (\sigma^*\gamma \cup \hat{\lambda}, -l_{\hat{\lambda}|\pi}),$$

by lemma 3.5, you get

$$[\sigma^*\gamma] \cup [(\hat{\lambda}, \hat{\mu})] = [(\sigma^*\gamma \cup \hat{\lambda}, \kappa^*\gamma \cup \hat{\mu})] = [(0, \kappa^*\gamma \cup \hat{\mu} + l_{\hat{\lambda}|\pi})]$$

in $H^3(\Delta, \pi; \mathbb{Z}/n)$. So, by the compatibility of $\exp_\pi$ and $\exp_\Delta$ and the perfect pairing

$$H^2(\Delta, \sigma^*A) \times H^1(\Delta, \pi; \sigma^*\hat{A}) \to H^3(\Delta, \pi; \sigma^*A),$$

the second sum reduces to

$$\sum_{[(\hat{\lambda},\hat{\mu})]} \exp_\pi\left(-l_{\hat{\lambda}|\pi} - \kappa^*\gamma \cup \hat{\mu}\right) = \sum_{[(\hat{\lambda},\hat{\mu})]} \exp_\Delta\left(-\sigma^*\gamma \cup (\hat{\lambda}, \hat{\mu})\right) = \mathbb{1}_{[\sigma^*\gamma]=0} \# H^1(\Delta, \pi; \sigma^*\hat{A}).$$

Therefore

$$\hat{x}_{\hat{\rho}}(\sigma) = \mathbb{1}_{[(\sigma^*\hat{\gamma}, \hat{\rho}^*\hat{a})]=0} \mathbb{1}_{[\sigma^*\gamma]=0} \# H^1(\Delta, \pi; \sigma^*\hat{A}) \exp_\pi(z_{\hat{\tau}|\pi} - \kappa^*\gamma \cup \hat{m} - c_{\hat{\rho}}).$$

By lemma 3.30,

$$\frac{\#A}{\#Z^1(\Delta, \sigma^*A)} = \frac{\#H^0(\Delta, \sigma^*A)}{\#H^1(\Delta, \sigma^*A)}$$

and by the perfect pairing,

$$\#H^1(\Delta, \pi; \sigma^*\hat{A}) = \#H^2(\Delta, \sigma^*A).$$

If $[(\sigma^*\hat{\gamma}, \hat{\rho}^*\hat{a})]$ or $[\sigma^*\gamma]$ is non-trivial, then $\hat{x}_{\hat{\rho}}(\sigma)$ and $x_{\hat{\rho}}(\sigma)$ both vanish. Otherwise

$$\frac{\hat{x}_{\hat{\rho}}(\sigma)}{x_{\hat{\rho}}(\sigma)} = \frac{\#H^1(\Delta, \pi; \sigma^*\hat{A}) \cdot \#A}{\#Z^1(\Delta, \sigma^*A)} \exp_\pi\left(z_{\hat{\tau}|\pi} - z_{\tau|\pi} - \kappa^*\gamma \cup \hat{m} + \tau^*_{|\pi}a \cup \hat{\rho}^*\hat{a}\right)$$

$$= \chi(\Delta, \sigma^*A) \exp_\pi\left(z_{\hat{\tau}|\pi} - z_{\tau|\pi} - \kappa^*\gamma \cup \hat{m} + \tau^*_{|\pi}a \cup (\hat{m}.\hat{\tau}_{|\pi})^*\hat{a}\right)$$

$$= \chi(\Delta, \sigma^*A) \exp_\pi\left(z_{\hat{\tau}|\pi} - z_{\tau|\pi} - \kappa^*\gamma \cup \hat{m} + \tau^*_{|\pi}a \cup (\hat{\tau}^*_{|\pi}\hat{a} - d\hat{m})\right)$$

$$= \chi(\Delta, \sigma^*A) \exp_\pi\left(z_{\hat{\tau}|\pi} - z_{\tau|\pi} + d(\tau^*_{|\pi}a \cup \hat{m}) + \tau^*_{|\pi}a \cup \hat{\tau}^*_{|\pi}\hat{a}\right)$$

$$= \chi(\Delta, \sigma^*A) \exp_\pi\left(z_{\hat{\tau}|\pi} - z_{\tau|\pi} + \tau^*_{|\pi}a \cup \hat{\tau}^*_{|\pi}\hat{a}\right)$$

$$= \chi(\Delta, \sigma^*A) \exp_\Delta(\hat{\tau}^*\hat{\omega} - \tau^*\omega + d(\tau^*a \cup \hat{\tau}^*\hat{a}), 0)$$
$$= \chi(\Delta, \sigma^*A) \exp_\Delta(\sigma^*\gamma \cup \hat{\tau}^*\hat{a} - \tau^*a \cup \sigma^*\hat{\gamma} + d(\tau^*a \cup \hat{\tau}^*\hat{a}), 0)$$
$$= \chi(\Delta, \sigma^*A) \exp_\Delta(0, 0) = \chi(\Delta, \sigma^*A) = \chi(\Delta, A)$$

finishes the proof. $\square$

## 4. Arithmetic Dijkgraaf-Witten theory

Let $G$ be a finite group and $n > 1$ a natural number. Let $F$ be a totally imaginary number field and $\zeta : \mathbb{Z}/n \xrightarrow{\sim} \mu_n$ an isomorphism over $F$. Set $X = \operatorname{spec} \mathcal{O}_F$ and choose a base point to define its étale fundamental group $\pi_1 X$. For the conventions regarding the cochain complex for group cohomology, see section 3.

**Definition 4.1.** The groupoid of fields over $X$ is
$$\mathcal{F}_X^G = \hom(\pi_1 X, G)//G$$
where $G$ acts by conjugation from the left.

The map $\varphi : X_{\text{ét}} \to X_{\text{fét}}$ induces a homomorphism
$$\varphi^* : H^3(\pi_1 X, \mathbb{Z}/n) \to H^3(X, \mathbb{Z}/n).$$

By Artin-Verdier duality [9, Theorem II.3.1], the invariant map induces an isomorphism
$$\operatorname{inv}_X : H^3(X, \mu_n) \xrightarrow{\sim} \operatorname{ext}_X^0(\mu_n, \mathbb{G}_m)^\vee \simeq \frac{1}{n}\mathbb{Z}/\mathbb{Z}$$
where $()^\vee$ is the Pontryagin dual. Consider the composition
$$\operatorname{tr}_X : H^3(\pi_1 X, \mathbb{Z}/n) \xrightarrow{\varphi^*} H^3(X, \mathbb{Z}/n) \xrightarrow{\zeta} H^3(X, \mu_n) \xrightarrow{\operatorname{inv}_X} \frac{1}{n}\mathbb{Z}/\mathbb{Z}.$$

Let $\omega \in Z^3(G, \mathbb{Z}/n)$ be an inhomogeneous, normalized cocycle.

**Definition 4.2.** Let $\tau : \pi_1 X \to G$ be an object of $\mathcal{F}_X^G$. The action of $\tau$ is
$$\mathcal{S}_X^\omega(\tau) = \operatorname{tr}_X(\tau^*\omega) \in \frac{1}{n}\mathbb{Z}/\mathbb{Z}.$$

**Lemma 4.3.** Let $\tau : \pi_1 X \to G$ be an object of $\mathcal{F}_X^G$. The action $\mathcal{S}_X^\omega(\tau)$ depends only on the equivalence class of $\tau$ and the cohomology class of $\omega$.

*Proof.* Since $\operatorname{tr}_X$ is defined on cohomology, the value $\mathcal{S}_X^\omega(\tau) = \operatorname{tr}_X(\varphi^*(\tau^*\omega))$ depends only on the cohomology class of $\omega$. Consider the chain homotopy
$$h_g : C^{\bullet+1}(G, \mathbb{Z}/n) \to C^\bullet(G, \mathbb{Z}/n)$$
between $\operatorname{ad}_g$ and $\operatorname{id}$ from lemma 3.8. Then
$$(g.\tau)^*\omega = \tau^*(\omega.g) = \tau^*(\omega + dh_g(\omega) + h_g(d\omega)) = \tau^*(\omega + dh_g(\omega)).$$

Since $dh_g(\omega)$ is a coboundary, $\mathcal{S}_X^\omega(\tau) = \mathcal{S}_X^\omega(g.\tau)$. $\square$

**Definition 4.4.** The arithmetic Dijkgraaf-Witten invariant of $X$ is
$$\mathcal{Z}^\omega(X) = \sum_{[\tau]} \exp(2\pi i \mathcal{S}_X^\omega(\tau)) \cdot \frac{1}{\#\operatorname{aut}(\tau)}$$

where the sum is over equivalence classes of objects $[\tau]$ of $\mathcal{F}_X^G$.

Define
$$\exp_X : H^3(\pi_1 X, \mathbb{Z}/n) \overset{\text{tr}_X}{\Rightarrow} \frac{1}{n}\mathbb{Z}/\mathbb{Z} \overset{\exp(2\pi i -)}{\hookrightarrow} \mathbb{C}^\times.$$

Consider the degree $\#G$ cover $\hom(\pi_1 X, G)$ of $\mathcal{F}_X^G$. Then the path integral equals
$$\mathcal{Z}^\omega(X) = \frac{1}{\#G} \sum_\tau \exp_X(\tau^* \omega)$$

where the sum is over $\tau \in \hom(\pi_1 X, G)$.

Consider an $n$-torsion sheaf $\mathcal{A}$ over $X$ and its dual $\hat{\mathcal{A}} = \hom(\mathcal{A}, \mathbb{Z}/n)$. Denote by $H^1(X, \mathcal{A})^\perp$ the orthogonal complement of $H^1(X, \mathcal{A})$ with respect to the cup product
$$H^1(X, \mathcal{A}) \times H^2(X, \hat{\mathcal{A}}) \overset{\cup}{\to} H^3(X, \mathbb{Z}/n).$$

Recall the construction of $A \rtimes_\gamma K$ from section 3.2.

**Theorem 4.5.** Let $n > 1$ be a natural number and
- $K$ a finite group, $A$ a finite, $n$-torsion $K$-module, and $\hat{A} = \hom(A, \mathbb{Z}/n)$
- $\gamma \in Z^2(K, A), \hat{\gamma} \in Z^2(K, \hat{A})$ inhomogeneous, normalized cocycles, and $e \in C^3(K, \mathbb{Z}/n)$ an inhomogeneous, normalized cochain such that $de = \gamma \cup \hat{\gamma}$
- $G = A \rtimes_\gamma K$ with projections $a$ and $k$ and $\hat{G} = \hat{A} \rtimes_{\hat{\gamma}} K$ with projections $\hat{a}$ and $\hat{k}$
- $\omega = k^* e + a \cup k^* \hat{\gamma} \in Z^3(G, \mathbb{Z}/n)$ and $\hat{\omega} = \hat{k}^* e + \hat{k}^* \gamma \cup \hat{a} \in Z^3(\hat{G}, \mathbb{Z}/n)$

Let $F$ be a totally imaginary number field and $\zeta : \mathbb{Z}/n \overset{\sim}{\to} \mu_n$ an isomorphism over $F$. Set $X = \spec \mathcal{O}_F$. For every homomorphism $\sigma : \pi_1 X \to K$, assume
$$[\sigma^* \gamma] \notin H^1(X, \sigma^* \hat{A})^\perp \setminus \{0\} \text{ and } [\sigma^* \hat{\gamma}] \notin H^1(X, \sigma^* A)^\perp \setminus \{0\}$$

as well as
$$\frac{\#H^1(X, \sigma^* A)}{\#H^0(X, \sigma^* A)} = \frac{\#H^1(X, \sigma^* \hat{A})}{\#H^0(X, \sigma^* \hat{A})}.$$

Then
$$\mathcal{Z}^\omega(X) = \mathcal{Z}^{\hat{\omega}}(X).$$

*Proof.* The path integral with respect to $\omega$ is
$$\mathcal{Z}^\omega(X) = \frac{1}{\#G} \sum_\tau \exp_X(\tau^*(k^* e + a \cup k^* \hat{\gamma})))$$

where the sum is over $\tau \in \hom(\pi_1 X, G)$. Since $\hom(\pi_1 X, G)$ is fibered over $\hom(\pi_1 X, K)$, the sum decomposes into
$$\mathcal{Z}^\omega(X) = \frac{1}{\#G} \sum_\sigma \sum_\tau \exp_X(\sigma^* e + \tau^* a \cup \sigma^* \hat{\gamma})$$

where the sum is over $\sigma \in \hom(\pi_1 X, K)$ and $\tau \in \{\tau : \pi_X \to G \mid k \circ \tau = \sigma\}$. For $\sigma : \pi_1 X \to K$, there is a homomorphism $\tau_\sigma : \pi_1 X \to G$ such that $k \circ \tau_\sigma = \sigma$ if and only if the class $[\sigma^* \gamma]$ vanishes in $H^2(\pi_1 X, \sigma^* A)$. If that is the case, choose such a $\tau_\sigma$. The choices are a $Z^1(\pi_1 X, \sigma^* A)$-torsor where $\lambda \in Z^1(\pi_1 X, \sigma^* A)$ acts by

$$\lambda \tau_\sigma = \lambda(a \circ \tau_\sigma, k \circ \tau_\sigma) = (a \circ \tau_\sigma - \lambda, \sigma).$$

So the path integral simplifies to

$$\mathcal{Z}^\omega(X) = \frac{1}{\#G} \sum_\sigma \mathbb{1}_{[\sigma^*\gamma]=0} \sum_\lambda \exp_X(\sigma^*e + (\lambda\tau_\sigma)^*a \cup \sigma^*\hat{\gamma})$$

$$= \frac{1}{\#G} \sum_\sigma \mathbb{1}_{[\sigma^*\gamma]=0} \sum_\lambda \exp_X(\sigma^*e + (\tau_\sigma^*a - \lambda) \cup \sigma^*\hat{\gamma})$$

$$= \frac{1}{\#G} \sum_\sigma \mathbb{1}_{[\sigma^*\gamma]=0} \exp(2\pi i \mathcal{S}^\omega(\tau_\sigma)) \sum_\lambda \exp_X(-\lambda \cup \sigma^*\hat{\gamma})$$

where the sum is over $\sigma \in \hom(\pi_1 X, K)$ and $\lambda \in Z^1(\pi_1 X, \sigma^*A)$. The pullback along the $\varphi : X_{\text{ét}} \to X_{\text{fét}}$ is an isomorphism in degree 1, injective in degree 2, and commutes with cup products. So by assumption, if $[\sigma^*\hat{\gamma}] \neq 0$, then $[\sigma^*\hat{\gamma}] \notin H^1(X, \sigma^*A)^\perp$. In this case

$$\exp_X(- \cup \sigma^*\hat{\gamma}) : Z^1(\pi_1 X, \sigma^*A) \to \mathbb{C}^\times$$

is a non-trivial character. So

$$\sum_\lambda \exp_X(-\lambda \cup \sigma^*\hat{\gamma}) = \mathbb{1}_{[\sigma^*\hat{\gamma}]=0} \# Z^1(\pi_1 X, \sigma^*A).$$

Therefore the path integral with respect to omega is

$$\mathcal{Z}^\omega(X) = \frac{1}{\#G} \sum_\sigma \mathbb{1}_{[\sigma^*\gamma]=0} \mathbb{1}_{[\sigma^*\hat{\gamma}]=0} \# Z^1(\pi_1 X, \sigma^*A) \exp_X(2\pi i \mathcal{S}^\omega(\tau_\sigma)).$$

Similarly, since $[\sigma^*\gamma] \notin H^1(X, \sigma^*\hat{A})^\perp \setminus \{0\}$, the path integral with respect to $\hat{\omega}$ is

$$\mathcal{Z}^{\hat{\omega}}(X) = \frac{1}{\#\hat{G}} \sum_\sigma \mathbb{1}_{[\sigma^*\hat{\gamma}]=0} \mathbb{1}_{[\sigma^*\gamma]=0} \# Z^1(\pi_1 X, \sigma^*\hat{A}) \exp_X(2\pi i \mathcal{S}^{\hat{\omega}}(\hat{\tau}_\sigma))$$

where the sum is over $\sigma \in \hom(\pi_1 X, K)$. By lemma 3.30,

$$\# Z^1(\pi_1 X, \sigma^*A) = \#A \frac{\# H^1(\pi_1 X, \sigma^*A)}{\# H^0(\pi_1 X, \sigma^*A)}$$

and

$$\# Z^1(\pi_1 X, \sigma^*\hat{A}) = \#A \frac{\# H^1(\pi_1 X, \sigma^*\hat{A})}{\# H^0(\pi_1 X, \sigma^*\hat{A})}.$$

Therefore, by assumption,

$$\# Z^1(\pi_1 X, \sigma^*A) = \# Z^1(\pi_1 X, \sigma^*\hat{A}).$$

Assume $[\sigma^*\gamma]$ and $[\sigma^*\hat{\gamma}]$ vanish. Then there is $\tau : \pi_1 X \to G$ and $\hat{\tau} : \pi_1 X \to \hat{G}$ such that $k \circ \tau = \hat{k} \circ \hat{\tau} = \sigma$. It remains to show $\mathcal{S}_X^\omega(\tau) = \mathcal{S}_X^{\hat{\omega}}(\hat{\tau})$. Consider $(\tau, \hat{\tau}) : \pi_1 X \to G \times_K \hat{G}$. Then, by implicitly composing with the pullbacks along the projections onto $G$ and $\hat{G}$, you get $\mathcal{S}^\omega(\tau) = \mathcal{S}^\omega(\tau, \hat{\tau})$ and $\mathcal{S}^{\hat{\omega}}(\hat{\tau}) = \mathcal{S}^{\hat{\omega}}(\tau, \hat{\tau})$. Over $G \times_K \hat{G}$,

$$d(a \cup \hat{a}) = da \cup \hat{a} - a \cup d\hat{a} = -k^*\gamma \cup \hat{a} + a \cup \hat{k}^*\hat{\gamma}.$$

Therefore the difference

$$\mathcal{S}^\omega(\tau) - \mathcal{S}^{\hat{\omega}}(\hat{\tau}) = \text{tr}_X(\sigma^*e + \tau^*a \cup \sigma^*\hat{\gamma}) - \text{tr}_X(\sigma^*e + \sigma^*\gamma \cup \hat{\tau}^*\hat{a})$$

$$= \mathrm{tr}_X(\tau^* a \cup \sigma^* \hat{\gamma} - \sigma^* \gamma \cup \hat{\tau}^* \hat{a})$$
$$= \mathcal{S}^{d(a \cup \hat{a})}(\tau, \hat{\tau}) = 0$$

vanishes. □

**Remark 4.6.** If $A$ is a trivial $K$-module, then, non-canonically, $A \simeq \hat{A}$ as $K$-modules. Therefore the second condition

$$\frac{\#H^1(X, \sigma^* A)}{\#H^0(X, \sigma^* A)} = \frac{\#H^1(X, \sigma^* \hat{A})}{\#H^0(X, \sigma^* \hat{A})}$$

of theorem 4.5 is satisfied for every $\sigma : \pi_1 X \to K$.

## 4.1. Boundaries

Let $n > 1$ be a natural number. Let $L_1, ..., L_r$ be non-archimedean local fields of characteristic zero and $\zeta_v : \mathbb{Z}/n \xrightarrow{\sim} \mu_n$ an isomorphism over $L_v$ for each $v$. Set $Y_v = \mathrm{spec}\, L_v$ and choose a base point to define its étale fundamental group $\pi_1 Y_v$. Consider

$$Y = \bigsqcup Y_v.$$

By local Tate duality [5, Theorem 7.2.6], the cohomological dimension of $Y_v$ is 2 and its cohomology groups are finite. Furthermore the local invariant map restricts to

$$\mathrm{inv}_v : H^2(Y_v, \mu_n) \xrightarrow{\sim} H^0(Y_v, \mu_n^D)^\vee \simeq \frac{1}{n}\mathbb{Z}/\mathbb{Z}$$

where $\mu_n^D = \mathscr{h}\mathit{om}(\mu_n, \mathbb{G}_m)$ is the Cartier dual. Define

$$\mathrm{tr}_{Y_v} : H^2(Y_v, \mathbb{Z}/n) \xrightarrow{\zeta_v} H^2(Y_v, \mu_n) \xrightarrow{\mathrm{inv}_v} \frac{1}{n}\mathbb{Z}/\mathbb{Z}.$$

Let $F$ be a totally imaginary number field and $\zeta : \mathbb{Z}/n \xrightarrow{\sim} \mu_n$ an isomorphism over $F$. Let $S$ be a finite set of finite places of $F$ containing all divisors of $n$. Suppose there is an isomorphism $\bigsqcup \mathrm{spec}\, L_v \simeq Y$ compatible with the trivializations $\zeta$ and $\zeta_v$. Set

$$X = \mathrm{spec}\, \mathcal{O}_F \setminus S$$

and choose a base point to define its étale fundamental group $\pi_1 X$. For each $v$, choose a path between the base point of $X$ and the base point of $Y_v$. By Artin-Verdier duality, the cohomological dimension of $X$ is 2 and its cohomology groups are finite. Furthermore the sum of the local invariant maps induce

$$\mathrm{inv}_X : H^3_c(X, \mu_n) \xrightarrow{\sim} \mathrm{ext}_X^0(\mu_n, \mathbb{G}_m)^\vee \simeq \frac{1}{n}\mathbb{Z}/\mathbb{Z}.$$

Define

$$\mathrm{tr}_X : H^3_c(X, \mathbb{Z}/n) \xrightarrow{\zeta} H^3_c(X, \mu_n) \xrightarrow{\mathrm{inv}_X} \frac{1}{n}\mathbb{Z}/\mathbb{Z}.$$

In particular

$$\begin{array}{ccc}
\bigoplus_{v \notin X} H^2(Y_v, \mathbb{Z}/n) & \longrightarrow & H^3_c(X, \mathbb{Z}/n) \\
{\scriptstyle \sum \mathrm{tr}_{Y_v}}\downarrow & & \downarrow{\scriptstyle \mathrm{tr}_X} \\
\frac{1}{n}\mathbb{Z}/\mathbb{Z} & =\!=\!=\!= & \frac{1}{n}\mathbb{Z}/\mathbb{Z}
\end{array}$$

commutes. Since $Y_v$ is the spectrum of a field, there is a canonical isomorphism

$$H^\bullet(\pi_1 Y_v, \mathbb{Z}/n) \xrightarrow{\sim} H^\bullet(Y_v, \mathbb{Z}/n).$$

Since $n$ is invertible on $X$, by [9, Proposition II.2.9], there is a canonical isomorphism

$$H^\bullet(\pi_1 X, \mathbb{Z}/n) \xrightarrow{\sim} H^\bullet(X, \mathbb{Z}/n).$$

By [7, Proposition 2.1], this isomorphism induces an isomorphism

$$H^\bullet(\pi_1 X, \pi_1 Y; \mathbb{Z}/n) \xrightarrow{\sim} H^\bullet_c(X, \mathbb{Z}/n),$$

which commutes with the traces. Therefore the pair of $\pi_1 X$ and $\pi_1 Y = (\pi_1 Y_v)$ satisfies the assumptions of section 3.1.

**Definition 4.7.** The arithmetic Dijkgraaf-Witten theory assigns to $Y$ the finite-dimensional complex Hilbert space

$$\mathcal{Z}^\omega(Y) = \mathcal{Z}^\omega(\pi_1 Y)$$

and to $X$ the vector

$$\mathcal{Z}^\omega(X) = \mathcal{Z}^\omega(\pi_1 X) \in \mathcal{Z}^\omega(Y).$$

**Remark 4.8.** Restricting to $n$ invertible is unsatisfying. If you look closer at the topological definition of Dijkgraaf-Witten theory, it is important to choose cochain complexes that are functorial in both $X$ and the coefficients. Moreover $\omega$ is chosen as a cocycle of $BG$ with values in $\mathbb{Q}/\mathbb{Z}$ twisted by the orientation sheaf. For arithmetic Dijkgraaf-Witten theory, a more natural approach would start with a cocycle on a stack $BG$ with values in $\mu_n$. Moreover theorem 4.9 requires a Poincaré-like duality theorem. This suggests working with fppf topology, where Artin-Verdier duality [7, Theorem 1.1] provides an appropriate duality even when $n$ is not invertible.

**Theorem 4.9.** Let $n > 1$ be a natural number and
- $K$ a finite group, $A$ a finite, $n$-torsion $K$-module, and $\hat{A} = \hom(A, \mathbb{Z}/n)$
- $\gamma \in Z^2(K, A), \hat{\gamma} \in Z^2(K, \hat{A})$ inhomogeneous, normalized cocycles, and $e \in C^3(K, \mathbb{Z}/n)$ an inhomogeneous, normalized cochain such that $de = \gamma \cup \hat{\gamma}$
- $G = A \rtimes_\gamma K$ with projections $a$ and $k$ and $\hat{G} = \hat{A} \rtimes_{\hat{\gamma}} K$ with projections $\hat{a}$ and $\hat{k}$
- $\omega = k^*e + a \cup k^*\hat{\gamma} \in Z^3(G, \mathbb{Z}/n)$ and $\hat{\omega} = \hat{k}^*e + \hat{k}^*\gamma \cup \hat{a} \in Z^3(\hat{G}, \mathbb{Z}/n)$

Let $L_1, ..., L_r$ be non-archimedean local fields of characteristic zero and $\zeta_v : \mathbb{Z}/n \xrightarrow{\sim} \mu_n$ an isomorphism over $L_v$ for each $v$. Set $Y_v = \spec L_v$ and choose a base point to define its étale fundamental group $\pi_1 Y_v$. Consider $Y = \bigsqcup Y_v$. Then there is an isomorphism

$$\Theta : \mathcal{Z}^\omega(Y) \xrightarrow{\sim} \mathcal{Z}^{\hat{\omega}}(Y).$$

Let $F$ be a totally imaginary number field and $\zeta : \mathbb{Z}/n \xrightarrow{\sim} \mu_n$ an isomorphism over $F$. Let $S$ be a finite set of finite places of $F$ containing all divisors of $n$. Suppose there is an isomorphism $\bigsqcup \spec L_v \simeq Y$ compatible with the trivializations $\zeta$ and $\zeta_v$. Set $X = \spec \mathcal{O}_F \setminus S$ and choose a base point to define its étale fundamental group $\pi_1 X$. For each $v$, choose a path between the base point of $X$ and the base point of $Y_v$. Then

$$\Theta(\mathcal{Z}^\omega(\Delta)) = \mathcal{Z}^{\hat{\omega}}(\Delta).$$

*Proof.* Suppose $\mathcal{A}$ is a finite, $n$-torsion $\pi_1 Y_v$-module. As before identify $H^\bullet(\pi_1 Y_v, \mathcal{A}) \simeq H^\bullet(Y_v, \mathcal{A})$. By [5, Theorem 7.3.1], the multiplicative Euler characteristic $\chi(\pi_1 Y_v, \mathcal{A})$ is independent of the $\pi_1 Y_v$-action. Since $\mathcal{A}$ is $n$-torsion, its Cartier dual is $\mathcal{A}^D = \hom(\mathcal{A}, \mu_n)$. So $H^i(\pi_1 Y_v, \mathcal{A}^D)$ is $n$-torsion. Therefore local Tate duality factors through

$$H^i(\pi_1 Y_v, \mathcal{A}) \times H^{2-i}(\pi_1 Y_v, \hom(\mathcal{A}, \mu_n)) \xrightarrow{\cup} H^2(\pi_1 Y_v, \mu_n).$$

Identify $\zeta_v : \mathbb{Z}/n \xrightarrow{\sim} \mu_n$ to see that

$$H^i(\pi_1 Y_v, \mathcal{A}) \times H^{2-i}(\pi_1 Y_v, \hat{\mathcal{A}}) \xrightarrow{\cup} H^2(\pi_1 Y_v, \mathbb{Z}/n)$$

is a perfect pairing of finite abelian groups for every integer $i$. Define

$$\Theta : \mathcal{Z}^\omega(Y) \xrightarrow{\sim} \mathcal{Z}^A(\pi_1 Y) \otimes \mathcal{Z}^\omega(Y) \xrightarrow{\sim} \mathcal{Z}^{\hat{\omega}}(Y)$$

where $\mathcal{Z}^A$ is defined in definition 3.20, the first map is multiplication by $\sqrt{\chi(\pi_1 Y, A)}$, and the second map is the isomorphism from theorem 3.23. Suppose $\mathcal{A}$ is a finite, $n$-torsion $\pi_1 X$-module. Again identify $\zeta : \mathbb{Z}/n \to \mu_n$. Since $n$ is invertible on $X$, by Artin-Verdier duality and [7, Proposition 6.2] (see [3, Corollary 6.11] for details), the cup product

$$H^i(\pi_1 X, \mathcal{A}) \times H^{3-i}(\pi_1 X, \pi_1 Y; \hat{\mathcal{A}}) \xrightarrow{\cup} H^3(\pi_1 X, \pi_1 Y; \mathbb{Z}/n)$$

is a perfect pairing of finite abelian groups for every integer $i$. By corollary 3.22 and theorem 3.23, the map satisfies $\Theta(\mathcal{Z}^\omega(X)) = \mathcal{Z}^{\hat{\omega}}(X)$. $\square$

### 4.2. Examples

Let $p_1, ..., p_r$ be prime numbers congruent to 1 modulo 4 such that their product is negative. Here a prime number is the choice of a generator of a non-zero prime ideal in $\mathbb{Z}$. For example choose $p_1 = 5$ and $p_2 = -7$. Denote their product by $d = p_1 \cdot ... \cdot p_r$ and consider the quadratic totally imaginary number field

$$F = \mathbb{Q}(\sqrt{d}).$$

**Lemma 4.10.** (Gauss genus theory) The class group of $F$ modulo 2 is isomorphic to

$$\mathrm{cl}(F)/2 \simeq \mathrm{gal}\big(F(\sqrt{p_1}, ..., \sqrt{p_r})/F\big) \simeq (\mathbb{Z}/2)^{r-1}.$$

*Proof.* Since $d$ is congruent to 1 modulo 4, the field $F$ is contained in $\mathbb{Q}(\mu_d)$. Consider the Galois action of the generator $\sigma$ of $\mathrm{gal}(F/\mathbb{Q})$ on $\mathrm{cl}(F)$. It is trivial on inert and ramified primes. If $\mathfrak{p}$ is split, then $\sigma(\mathfrak{p}) \cdot \mathfrak{p} = (p)$ is trivial in the class group. So the class $[\sigma(\mathfrak{p})]$ is inverse to $[\mathfrak{p}]$. Hence $\mathrm{gal}(F/\mathbb{Q})$ acts trivially on $\mathrm{cl}(F)/2$. Denote by $M/F$ the field extension corresponding to $\mathrm{cl}(F)/2$. Since it is the product over all unramified 2-extensions of $F$, it is Galois over $\mathbb{Q}$. So

$$0 \to \mathrm{gal}(M/F) \to \mathrm{gal}(M/\mathbb{Q}) \to \mathrm{gal}(F/\mathbb{Q}) \to 0$$

is a central extension of a cyclic group by an abelian extension. Therefore $\mathrm{gal}(M/\mathbb{Q})$ is abelian. Since the ramification of $M$ over $\mathbb{Q}$ is tame and has support in the divisors of $d$, the field $M$ is contained in $\mathbb{Q}(\mu_d)$. Identify $\mathrm{gal}(\mathbb{Q}(\mu_d)/\mathbb{Q})$ with $\prod (\mathbb{Z}/p_i)^\times$ to get

$$\begin{array}{ccccccccc}
0 & \longrightarrow & K & \longrightarrow & \prod (\mathbb{Z}/p_i)^\times & \longrightarrow & \{\pm 1\} & \longrightarrow & 0 \\
& & \downarrow & & \downarrow & & \downarrow \wr & & \\
0 & \longrightarrow & \mathrm{gal}(M/F) & \longrightarrow & \mathrm{gal}(M/\mathbb{Q}) & \longrightarrow & \mathrm{gal}(F/\mathbb{Q}) & \longrightarrow & 0
\end{array}$$

where $K$ is the kernel of the product of the Legendre symbols. Since at least one of the primes $p_i$ satisfies $|p_i| \equiv 3 \bmod 4$, the top row splits. Therefore the first exact sequence modulo 2

$$0 \to K/2 \to \prod_{i=1}^r \{\pm 1\} \to \{\pm 1\} \to 0$$

is exact. Since $\text{gal}(M/F) \simeq \text{cl}(F)/2$, the first vertical homomorphism factors through $K/2$. The factors of $\prod\{\pm 1\}$ correspond to the field extensions $\mathbb{Q}(\sqrt{p_i})/\mathbb{Q}$. By Abhyankar's lemma, those are contained in $M$, so there is a surjective homomorphism $\text{gal}(M/\mathbb{Q}) \to \prod\{\pm 1\}$. This already implies that the bottom row can be identified with this exact sequence. Therefore $\text{gal}(M/F)$ is canonically isomorphic to $K/2$. □

**Lemma 4.11.** Let $\mathfrak{p}_i$ be the prime in $F$ above $p_i$. Then the 2-torsion $\text{cl}(F)[2]$ of the class group is generated by the classes $[\mathfrak{p}_1], \ldots, [\mathfrak{p}_r]$ with a single relation $[\mathfrak{p}_1] + \ldots + [\mathfrak{p}_r] = 0$.

*Proof.* The following proof is from a short note by Guo, which is based on [10]. Take a class $[\mathfrak{a}] \in \text{cl}(F)[2]$ and write $\mathfrak{a}^2 = (a)$. Then

$$\text{norm}_{F/\mathbb{Q}}(a) = \text{norm}_{F/\mathbb{Q}}(\mathfrak{a}^2) = b^2$$

for some $b \in \mathbb{Q}$. Since $\text{norm}_{F/\mathbb{Q}}(a/b) = 1$, by Hilbert 90, there is $c \in F^\times$ such that

$$\frac{a}{b} = \frac{c}{\sigma(c)}.$$

So

$$a = \frac{b \cdot \text{norm}_{F/\mathbb{Q}}(c)}{\sigma(c)^2}.$$

Define $\mathfrak{a}' = \mathfrak{a} \cdot \sigma(c)$. Then the classes $[\mathfrak{a}]$ and $[\mathfrak{a}']$ coincide and

$$\mathfrak{a}'^2 = (b \cdot \text{norm}_{F/\mathbb{Q}}(c)) = (a')$$

is a rational number. Since $\mathfrak{a}'$ can be multiplied by any rational number without changing its ideal class, assume $a'$ is square free. Since $\text{div}(\mathfrak{a}') = \text{div}(a')/2$ and $a'$ is square free, the ideal $\mathfrak{a}'$ has support in the ramified primes. So $[\mathfrak{p}_1], \ldots, [\mathfrak{p}_r]$ generate the 2-torsion of the class group. By lemma 4.10, the rank of $\text{cl}(F)[2]$ is $r - 1$. Therefore the relation

$$\mathfrak{p}_1 \cdot \ldots \cdot \mathfrak{p}_r = \left(\sqrt{d}\right)$$

is enough. □

**Lemma 4.12.** Let

$$k = \begin{cases} 0 & \text{if } n \text{ odd} \\ n/2 & \text{if } n \text{ even} \end{cases}$$

Then

$$H^\bullet(\mathbb{Z}/n, \mathbb{Z}/n) \simeq \mathbb{Z}/n[x, y]/(x^2 - ky)$$

where $x$ is of degree 1 and $y$ is of degree 2.

*Proof.* See [11, Example 3.41]. □

Since 2 is prime, by the Künneth formula,

$$H^\bullet(\mathbb{Z}/2^s, \mathbb{Z}/2) = \bigotimes_{i=1}^{s} H^\bullet(\mathbb{Z}/2, \mathbb{Z}/2) = \bigotimes_{i=1}^{s} \mathbb{Z}/2[x_i]$$

where $x_i$ is the projection onto the $i$-th factor.

Set $X = \text{spec } \mathcal{O}_F$.

**Example 4.13.** Consider the quaternion group $G = Q_8$ and the trivial cocycle $\omega = 0 \in Z^3(Q_8, \mathbb{Z}/2)$. The central extension

$$0 \to \mathbb{Z}/2 \to Q_8 \to (\mathbb{Z}/2)^2 \to 0$$

is represented by the cocycle

$$\gamma = x \cup y + x \cup x + y \cup y \in Z^2\big((\mathbb{Z}/2)^2, \mathbb{Z}/2\big)$$

where $x, y : (\mathbb{Z}/2)^2 \to \mathbb{Z}/2$ are the projections. The arithmetic Dijkgraaf-Witten invariant is

$$\mathcal{Z}^\omega(X) = \frac{1}{8} \# \hom(\pi_1 X, Q_8).$$

For the dual, consider $\hat{G} = (\mathbb{Z}/2)^3$ as the trivial extension

$$0 \to \mathbb{Z}/2 \to (\mathbb{Z}/2)^3 \to (\mathbb{Z}/2)^2 \to 0$$

and

$$\hat{\omega} = \gamma \cup z = x \cup y \cup z + x \cup x \cup z + y \cup y \cup z \in Z^3\big((\mathbb{Z}/2)^3, \mathbb{Z}/2\big)$$

where $x, y, z : (\mathbb{Z}/2)^3 \to \mathbb{Z}/2$ are the projections. View $H^3(X, \mathbb{Z}/2) = \mu_2^\vee$. Then the path integral is

$$\mathcal{Z}^{\hat{\omega}}(X) = \frac{1}{8} \sum_{\hat{\tau}} \exp_X(\hat{\tau}^* \hat{\omega})$$

where the sum is over homomorphisms $\hat{\tau} : \text{cl}(F)/2 \to (\mathbb{Z}/2)^3$. Choose $e = 0$ to see this in the context of theorem 4.5. So, if $[\sigma^* \gamma]$ either vanishes or is not contained in $H^1(X, \mathbb{Z}/2)$ for every homomorphism $\sigma : \pi_1 X \to (\mathbb{Z}/2)^2$, duality $\mathcal{Z}^\omega(X) = \mathcal{Z}^{\hat{\omega}}(X)$ holds. The remaining part is about understanding this condition better.

The first step is to better understand $H^2(X, \mathbb{Z}/2)$. By [12, Corollary 2.15], elements of $H^2(X, \mathbb{Z}/2)$ should be seen as functions on the dual

$$H^2(X, \mathbb{Z}/2)^\vee = \{(\mathfrak{a}, a) \in \text{cl}(F) \oplus F^\times \mid 2\mathfrak{a} + (a) = 0\} / \{((a), a^{-2}) \mid a \in F^\times\}.$$

Since $F$ is a quadratic, totally imaginary number field, by Dirichlet's unit theorem, $\mathcal{O}_F^\times/2 \simeq \{\pm 1\}$. So there is an exact sequence

$$0 \to \{\pm 1\} \to H^2(X, \mathbb{Z}/2)^\vee \to \text{cl}(F)[2] \to 0.$$

Next examine $H^1(X, \mathbb{Z}/2)$. Denote by $x_i \in H^1(X, \mathbb{Z}/2)$ the class corresponding to the field extension $F(\sqrt{p_i})$ over $F$. By lemma 4.10, the classes $x_i$ for $i = 1, ..., r-1$ form a basis of $H^1(X, \mathbb{Z}/2)$. Denote by $\mathfrak{p}_i$ the prime in $F$ above $p_i$. Consider the isomorphism

$$t : H^1(X, \mathbb{Z}/2) \xrightarrow{\sim} \text{cl}(F)[2]$$
$$x_i \mapsto \mathfrak{p}_i$$

and the section

$$s : \text{cl}(F)[2] \hookrightarrow H^2(X, \mathbb{Z}/2)^\vee$$
$$\mathfrak{p}_i \mapsto (\mathfrak{p}_i, p_i^{-1}).$$

**Lemma 4.14.** The adjoint map of the trace composed with the cup product

$$H^1(X, \mathbb{Z}/2) \times H^2(X, \mathbb{Z}/2) \overset{\cup}{\to} H^3(X, \mathbb{Z}/2) \overset{\text{tr}_X}{\to} \frac{1}{2}\mathbb{Z}/\mathbb{Z}$$

is the composition

$$s \circ t : H^1(X, \mathbb{Z}/2) \to H^2(X, \mathbb{Z}/2)^{\vee}.$$

*Proof.* Take $x_i \in H^1(X, \mathbb{Z}/2)$ as before. Denote by $L = F(\sqrt{p_i})$ the corresponding field extension and by $\tau$ the generator of the Galois group gal$(L/F)$. Then

$$-1 = \tau(\sqrt{p_i})/\sqrt{p_i}$$
$$\mathfrak{p}_i \mathcal{O}_L = (\sqrt{p_i}).$$

View $H^3(X, \mathbb{Z}/2) = \mu_2^{\vee}$. Then, by [12, Proposition 1.2],

$$(z \cup u)(-1) = u(\mathfrak{p}_i, p_i^{-1}),$$

which shows the claim. $\square$

Now examine double cup products. Take $x, y \in H^1(X, \mathbb{Z}/2)$ and denote by $M, N$ the corresponding field extensions. Use [12, Proposition 1.1] to calculate $(x \cup y)(\mathfrak{a}, a)$ as follows. Denote by $\sigma$ the generator of gal$(M/F)$. The computationally heavy part is to find an element $b \in M$ and an ideal $\mathfrak{b}$ such that

$$\mathfrak{a}\mathcal{O}_M = (1 - \sigma)\mathfrak{b} + (b)$$
$$\text{norm}_{M/F}(b) = a^{-1} \bmod \text{norm}_{M/F}(\mathcal{O}_M^{\times})$$

Then

$$(x \cup y)(\mathfrak{a}, a) = \text{artin}_{N/F}\big(\text{norm}_{M/F}(\mathfrak{b}) + \mathfrak{a}\big).$$

Here the Artin symbol counts the inert prime divisors modulo 2. If $x = y$, since

$$\text{norm}_{M/F} \text{cl}\, M = \ker \text{artin}_{M/F},$$

you can skip the hard step and get

$$(x \cup x)(\mathfrak{a}, a) = \text{artin}_{N/F}(\mathfrak{a}),$$

which depends only on the class $\mathfrak{a} \in \text{cl}(F)[2]$.

**Remark 4.15.** Applying these methods to example 4.13 using the Hecke library [13] results in (see https://github.com/jaroeichler/thesis):

| Primes | $\mathcal{Z}^{\omega}(X)$ | $\mathcal{Z}^{\hat{\omega}}(X)$ | Linking form |
|---|---|---|---|
| $-11, -83, -107, -139, -191$ | 8 | 8 | Non-symmetric |
| $29, -31, -43, -47, 101$ | 8 | 20 | Non-symmetric |
| $-11, -59, -107$ | 0.5 | 3.5 | Non-symmetric |
| $5, 193, -439$ | 0.5 | 0.5 | Symmetric |

Calculating further examples, there are two observations:
1. If the linking form is symmetric, then $\mathcal{Z}^{\omega}(X) = \mathcal{Z}^{\hat{\omega}}(X)$
2. If $|p_1|, ..., |p_{r-1}| \equiv 1 \bmod 4$, then the linking form is symmetric

Similarly consider the dihedral group with 8 elements $G = D_4$, which is represented by $\gamma = x \cup y + y \cup y$ and $\omega = 0$. Compare this with $\hat{G} = (\mathbb{Z}/2^3)$ and $\hat{\omega} = \gamma \cup z$. Here the calculations suggest $\mathcal{Z}^\omega(X) = \mathcal{Z}^{\hat\omega}(X)$ always holds.

The remaining part shows the first two observations from remark 4.15.

**Lemma 4.16.** The Bockstein is a derivation. So

$$\beta(x \cup y) = \beta(x) \cup y + (-1)^{\deg(x)} x \cup \beta(y).$$

*Proof.* See [14, Lemma 2.1]. □

**Lemma 4.17.** Let $x \in H^1(X, \mathbb{Z}/2)$. Then

$$\beta(x) = x \cup x.$$

*Proof.* Consider $x \in H^1(X, \mathbb{Z}/2) \simeq H^1(\pi_1 X, \mathbb{Z}/2)$ as a homomorphism $\chi : \pi_1 X \to \mathbb{Z}/2$. Since

$$\begin{array}{ccc} H^1(X, \mathbb{Z}/2) & \xrightarrow{\beta} & H^2(X, \mathbb{Z}/2) \\ \uparrow\wr & & \uparrow\wr \\ H^1(\pi_1 X, \mathbb{Z}/2) & \xrightarrow{\beta} & H^2(\pi_1 X, \mathbb{Z}/2) \\ \uparrow \chi^* & & \uparrow \chi^* \\ H^1(\mathbb{Z}/2, \mathbb{Z}/2) & \xrightarrow{\beta} & H^2(\mathbb{Z}/2, \mathbb{Z}/2) \end{array}$$

commutes and sends the non-trivial element $\alpha \in H^1(\mathbb{Z}/2, \mathbb{Z}/2)$ to $x$, it is enough to show that $\beta(\alpha) = \alpha \cup \alpha$. The group $H^i(\mathbb{Z}/2, \mathbb{Z})$ vanishes for $i$ odd and is isomorphic to $\mathbb{Z}/2$ for $i > 0$ even. So the exact sequence

$$0 \to \mathbb{Z} \xrightarrow{2} \mathbb{Z} \to \mathbb{Z}/2 \to 0$$

shows that $\beta : H^1(\mathbb{Z}/2, \mathbb{Z}/2) \xrightarrow{\sim} H^2(\mathbb{Z}/2, \mathbb{Z}/2)$ is an isomorphism. By [11, Example 3.9], the cup product $\alpha \cup \alpha$ is non-zero. This already implies $\beta(\alpha) = \alpha \cup \alpha$. □

In analogy to the linking form in topology, consider the following definition.

**Definition 4.18.** The linking form is

$$L_X : H^1(X, \mathbb{Z}/2) \times H^1(X, \mathbb{Z}/2) \to H^3(X, \mathbb{Z}/2) \xrightarrow{\sim} H^3(X, \mu_2) \xrightarrow{\text{inv}_X} \frac{1}{2}\mathbb{Z}/\mathbb{Z}$$

$$x, y \mapsto x \cup \beta(y)$$

By lemma 4.17, the linking form for $X$ with coefficients in $\mathbb{Z}/2$ can be calculated as

$$L_X(x, y) = x \cup y \cup y.$$

**Lemma 4.19.** The linking form $L_X$ is symmetric if and only if

$$(x \cup y)(\mathcal{O}_F, -1) = 0$$

for all $x, y \in H^1(X, \mathbb{Z}/2)$.

*Proof.* Denote by $\beta$ the Bockstein of the sequence

$$0 \to \mathbb{Z}/2 \to \mathbb{Z}/4 \to \mathbb{Z}/2 \to 0.$$

By lemma 4.16, $\beta$ is a derivation and by lemma 4.17, $\beta(x) = x \cup x$. Identify $H^3(X, \mathbb{Z}/2) = \{\pm 1\}^\vee$. Then

$$\beta : H^2(X, \mathbb{Z}/2) \to H^3(X, \mathbb{Z}/2)$$

is dual to

$$\{\pm 1\} \to H^2(X, \mathbb{Z}/2)^\vee$$
$$-1 \mapsto (\mathcal{O}_F, -1).$$

Therefore

$$(x \cup y)(\mathcal{O}_F, -1) = \beta(x \cup y)(-1) = (\beta(x) \cup y - x \cup \beta(y))(-1)$$
$$= (y \cup \beta(x) - x \cup \beta(y))(-1) = L_X(y, x) - L_X(x, y)$$

is exactly the symmetry condition. $\square$

**Lemma 4.20.** The cup product

$$H^1(X, \mathbb{Z}/2) \times H^2(X, \mathbb{Z}/2) \overset{\cup}{\to} H^3(X, \mathbb{Z}/2)$$

is left non-degenerate and the orthogonal complement $H^1(X, \mathbb{Z}/2)^\perp$ is a 1-dimensional $\mathbb{F}_2$-vector space.

*Proof.* The pairing is left non-degenerate if and only if the adjoint map $H^1(X, \mathbb{Z}/2) \to H^2(X, \mathbb{Z}/2)^\vee$ is injective. By lemma 4.14, this map is the composition

$$H^1(X, \mathbb{Z}/2) \overset{t}{\hookrightarrow} \mathrm{cl}(F)[2] \overset{s}{\hookrightarrow} H^2(X, \mathbb{Z}/2)^\vee,$$

which is injective. The exact sequence

$$0 \to \{\pm 1\} \to H^2(X, \mathbb{Z}/2)^\vee \to \mathrm{cl}(F)[2] \to 0$$

shows

$$\dim_{\mathbb{F}_2} H^2(X, \mathbb{Z}/2) = \dim_{\mathbb{F}_2} H^1(X, \mathbb{Z}/2) + 1.$$

So $H^1(X, \mathbb{Z}/2)^\perp$ is 1-dimensional. $\square$

**Lemma 4.21.** Let $u \in H^2(X, \mathbb{Z}/2)$ such that $u(\mathcal{O}_F, -1) = 0$. Then $u \notin H^1(X, \mathbb{Z}/2)^\perp \setminus \{0\}$.

*Proof.* The section $s$ of

$$0 \to \mathcal{O}_F^\times/2 \to H^2(X, \mathbb{Z}/2)^\vee \to \mathrm{cl}(F)[2] \to 0$$

induces

$$H^2(X, \mathbb{Z}/2) \simeq \left(\mathcal{O}_F^\times/2\right)^\vee \oplus \mathrm{cl}(F)[2]^\vee.$$

Then lemma 4.14 shows the claim. $\square$

The following corollary shows the first observation of remark 4.15.

**Corollary 4.22.** Let $\omega$ and $\hat{\omega}$ be as in example 4.13. Assume the linking form $L_X$ is symmetric. Then

$$\mathcal{Z}^\omega(X) = \mathcal{Z}^{\hat{\omega}}(X).$$

*Proof.* Recall $\gamma = x \cup y + x \cup x + y \cup y$ from example 4.13. Since $L_X$ is symmetric and pullbacks commute with cup products, by lemma 4.19, $\sigma^* \gamma(\mathcal{O}_F, -1) = 0$ for every $\sigma : \pi_1 X \to (\mathbb{Z}/2)^2$. So, by lemma 4.21, the assumptions of theorem 4.5 are satisfied, which shows the claim. $\square$

**Lemma 4.23.** Let $L$ be a number field, $M, N$ quadratic field extensions of $L$ such that $M \cap N = L$, and $\mathfrak{n}$ an ideal of $\mathcal{O}_N$ that is unramified over $L$. Then
$$\text{norm}_{N/L}(\mathfrak{n})\mathcal{O}_M = \text{norm}_{MN/M}(\mathfrak{n}\mathcal{O}_{MN}).$$

*Proof.* It is enough to show the claim for $\mathfrak{n}$ prime. Denote the prime in $L$ below $\mathfrak{n}$ by $\mathfrak{p}$. If $M$ is ramified over $\mathfrak{p}$, then $MN$ is ramified over $\mathfrak{n}$. Denote the prime in $M$ above $\mathfrak{p}$ by $\mathfrak{m}$. If $\mathfrak{p}$ is inert in $N$ and $\mathfrak{q}$ is the prime over $\mathfrak{n}$ in $MN$, then
$$\text{norm}_{N/L}(\mathfrak{n})\mathcal{O}_M = 2\mathfrak{p}\mathcal{O}_M = 4\mathfrak{m} = \text{norm}_{MN/M}(2\mathfrak{q}) = \text{norm}_{MN/M}(\mathfrak{n}\mathcal{O}_{MN}).$$

The picture to have in mind is the square

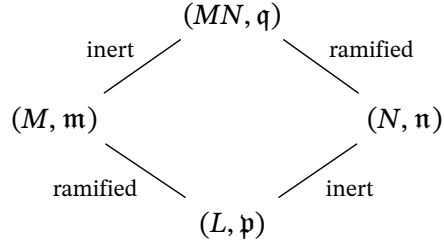

A similar calculation shows the other cases, except for when $\mathfrak{p}$ splits in $M$ and $N$ but is not totally split in $MN$, and $\mathfrak{p}$ is inert in $M$ and $N$ but not in $MN$. So it is enough to show that if $\mathfrak{p}$ splits in $N$ and $\mathfrak{m}$ is a prime in $M$ above $\mathfrak{p}$, then $\mathfrak{m}$ splits in $MN$. Let $\mathfrak{n}'$ be the other prime in $N$ above $\mathfrak{p}$. Then
$$L_\mathfrak{p} \otimes_L N = N_\mathfrak{n} \oplus N_{\mathfrak{n}'}.$$

Since $M$ and $N$ have trivial intersection,
$$M_\mathfrak{m} \otimes_M MN = M_\mathfrak{m} \otimes_L N = M_\mathfrak{m} \otimes_{L_\mathfrak{p}} L_\mathfrak{p} \otimes_L N = M_\mathfrak{m} \otimes_{L_\mathfrak{p}} N_k \oplus M_\mathfrak{m} \otimes_{L_\mathfrak{p}} N_{k'}$$

shows the claim. $\square$

The following lemma shows the second observation of remark 4.15.

**Lemma 4.24.** Let $K$ be a finite group and $A$ a finite, trivial, 2-torsion $K$-module. Assume
$$|p_1|, ..., |p_{r-1}| \equiv 1 \bmod 4.$$

Then the linking form $L_X$ is symmetric.

*Proof.* By lemma 4.10, any unramified degree 2 extension of $F$ is of the form $F(\sqrt{p_I})$ where $p_I$ is the product of $p_i$ for $i \in I$ and $I \subseteq \{1, ..., r-1\}$. Consider
$$M = \mathbb{Q}(\sqrt{p_I}).$$

By the Hasse norm theorem, $-1$ is a norm of $M/\mathbb{Q}$ if and only if it is locally a norm at every place. At unramified places $-1$ is always a norm. So it is enough to examine the primes $p_i$ for $i \in I$. If $-1$ is a norm in the residue field extension, it lifts to a norm in the extension of local fields. In the residue field of a ramified prime, the norm is squaring. So it is enough to ensure that $-1$ is a square modulo $p_i$ for every $i \in I$. Since $|p_i| \equiv 1 \bmod 4$, the Legendre symbol of $-1$ in $p_i$ is 1, so $-1$ is a square. Take an element $u \in M$ such that $\text{norm}_{M/\mathbb{Q}}(u) = -1$. Then
$$\text{norm}_{FM/F}(u) = -1.$$

Since $-1$ is a unit, $u$ has support in the split primes. Consider the divisor of zeroes $\mathfrak{a}$ of $u$, then $\mathfrak{a} - \sigma_M(\mathfrak{a}) = \text{div}(u)$ and

$$\mathfrak{a}\mathcal{O}_{FM} - \sigma_{FM}(\mathfrak{a}\mathcal{O}_{FM}) = \text{div}(u).$$

Since $\mathfrak{a}$ is unramified over $\mathbb{Q}$, by lemma 4.23,

$$\text{norm}_{M/\mathbb{Q}}(\mathfrak{a})\mathcal{O}_F = \text{norm}_{FM/F}(\mathfrak{a}\mathcal{O}_{FM}).$$

Take the cohomology class $x \in H^1(X, \mathbb{Z}/2)$ corresponding to $FM$ and $y \in H^1(X, \mathbb{Z}/2)$ corresponding to a field extension $N/F$. Since $\text{norm}_{M/\mathbb{Q}}(\mathfrak{a})$ is a principal ideal,

$$(x \cup y)(\mathcal{O}_F, -1) = \text{artin}_{FN/F}\big(\text{norm}_{FM/F}(\mathfrak{a}\mathcal{O}_{FM})\big) = \text{artin}_{FN/F}\big(\text{norm}_{M/\mathbb{Q}}(\mathfrak{a})\mathcal{O}_F\big) = 0.$$

So, by lemma 4.19, the linking form is symmetric. $\square$